\magnification=\magstep1
\input amstex
\documentstyle{amsppt}
\hoffset=.25truein
\hsize=6truein
\vsize=8.75truein

\topmatter
\centerline{\bf ON THE TRACE OF HECKE OPERATORS FOR}
\title
 MAASS FORMS FOR CONGRUENCE SUBGROUPS
\endtitle
\abstract
  Let $E_\lambda$ be a Hilbert space, whose elements are functions
spanned by the eigenfunctions of the Laplace-Beltrami operator
associated with an eigenvalue $\lambda>0$.  The norm of elements
in this space is given by the Petersson inner product.
In this paper, the trace of Hecke operators $T_n$ acting on the
space $E_\lambda$ is computed for congruence subgroups $\Gamma_0(N)$
of square free level, which may be considered as the analogue of
the Eichler-Selberg trace formula \cite{11} for non-holomorphic
cusp forms of weight zero.
\endabstract
\author
J. B. Conrey and Xian-Jin Li
\endauthor
\address
American Institute of Mathematics, 360 Portage Avenue, Palo Alto, CA 94306
\endaddress
\thanks
Research of both authors supported by the American Institute of Mathematics.
\endthanks
\email
conrey\@aimath.org, \,\, xianjin\@math.Stanford.EDU
\endemail
\endtopmatter
\address
Current address for Xian-Jin Li: Department of Mathematics,
Brigham Young University, Provo, Utah 84602 USA
\endaddress
\document

\heading
1.   Introduction
\endheading

   Let $N$ be a positive integer greater than one.  Denote by
$\Gamma_0(N)$ the Hecke congruence subgroup of level $N$.
The Laplace-Beltrami  operator $\Delta$ on the upper half-plane
$\Cal H$ is given by
$$\Delta= -y^2\left(\frac{\partial ^2}{\partial x^2}+\frac{\partial
^2}{\partial y^2} \right).$$
Let $D$ be the fundamental domain of $\Gamma_0( N)$.   Eigenfunctions
of the discrete spectrum of $\Delta$ are nonzero
real-analytic solutions of the equation
$$\Delta\psi=\lambda\psi$$
such that $\psi(\gamma z)=\psi(z)$ for all $\gamma$ in $\Gamma_0( N)$ and
such that
$$\int_D |\psi(z)|^2dz<\infty$$
where $dz$ represents the Poincar\'e measure of the upper half-plane.

  Let $\frak a$ be a cusp of $\Gamma_0(N)$.  Its stabilizer is denoted
by $\Gamma_{\frak a}$.  An element $\sigma_{\frak a}\in PSL(2, \Bbb R)$
exists such that $\sigma_{\frak a}\infty=\frak a$ and
$\sigma_{\frak a}^{-1}\Gamma_{\frak a}\sigma_{\frak a}=\Gamma_\infty$.
Let $f$ be a $\Gamma_0(N)$-invariant function.  If $\frak a$ is a cusp
of $\Gamma_0(N)$, then $f(\sigma_{\frak a}z)$ is $\Gamma_\infty$-invariant,
and hence it admits a formal Fourier expansion
$$f(\sigma_{\frak a}z)=\sum_{n\in\Bbb Z}c_{\frak a n}(y)e^{2n\pi ix}.$$
A $\Gamma_0(N)$-invariant function is said to be a Maass cusp form if it is
square-integrable and is an eigenvector of $\Delta$, such that the Fourier
coefficient $c_{\frak a 0}(y)=0$ for every cusp $\frak a$ of $\Gamma_0(N)$.
If $\psi$ is a cusp form associated with a positive discrete eigenvalue
$\lambda$, then it has the Fourier expansion \cite {5}
$$\psi(\sigma_{\frak a}z)=\sqrt y\sum_{m\neq 0}\rho_{\frak a}(m)
K_{i\kappa}(2\pi |m|y)e^{2m\pi ix},$$
where $\kappa=\sqrt{\lambda-1/4}$ and $K_\nu(y)$ is given by
the formula $\S$6.32, \cite {18}
$$\aligned K_\nu(y)&=\frac{2^\nu\Gamma(\nu+\frac 12)}{y^\nu\sqrt\pi}\int_0^\infty
\frac{\cos (yt)}{(1+t^2)^{\nu+\frac 12}}dt\\
&=\frac 12 \int_0^\infty \exp\left(-\frac y2(t+\frac 1t)
\right)\frac{dt}{t^{\nu+1}}.\endaligned \tag 1.1$$
The complex numbers $\rho_{\frak a}(m)$, $m(\neq 0)\in \Bbb Z$, are called
the Fourier coefficients of $\psi$ around the cusp $\frak a$.

   The Hecke operators $T_n$, $n=1, 2, \cdots$, $(n, N)=1$, which act
in the space of automorphic functions with respect to $\Gamma_0( N)$,
are defined by
$$\left(T_nf\right)(z)=\frac 1{\sqrt n}\sum_{ad=n, \,
0\leq b<d}f\left(\frac {az+b}d\right).$$
An orthonormal system of eigenfunctions of $\Delta$ exists \cite {5} such
that each of them is an eigenfunction of all the Hecke operators.  We call
these eigenfunctions Maass-Hecke eigenfunctions.
Let $\lambda_j$, $j=1, 2, \cdots$, be an enumeration in increasing order
of all positive discrete eigenvalues of $\Delta$ for $\Gamma_0(N)$
with an eigenvalue of multiplicity $m$ appearing $m$ times, and let
$\kappa_j=\sqrt{\lambda_j-1/4}$.   If $\psi_j(z)$ is a Maass-Hecke
 eigenfunction of $\Delta$ associated with the $j$th eigenvalue
$\lambda_j $, then
$$\left(T_n\psi_j\right)(z)=\tau_j(n)\psi_j(z)$$
where $\rho_{j\infty}(m)=\rho_{j\infty}(d)\tau_j(n)$ if $m=dn$
with $n\geqslant 1$, $(n, dN)=1$.

  Let $E_\lambda$ be a Hilbert space
of functions spanned by the eigenfunctions of $\Delta$ with a positive
eigenvalue $\lambda$.  The inner product of the space is given by
$$\langle F(z), G(z)\rangle=\int_D F(z)\bar G(z)dz. \tag 1.2$$
The analogue for Maass forms of Eichler-Selberg's trace formula
\cite {11}, p.85 for modular forms is obtained for the full
modular group in \cite {8}.
In this paper, the trace $\text{tr} T_n$ of Hecke operators acting
on the space $E_\lambda$ is computed for congruence subgroups
$\Gamma_0(N)$.  Some of this computation is implicit in
Hejhal \cite{3}.

    Denote by $h_d$ the class number of
indefinite rational quadratic forms with discriminant d.  Define
$$\epsilon_d=\frac{v_0+u_0\sqrt d}2 \tag 1.3$$
where the pair $(v_0, u_0)$ is the fundamental solution \cite {9} of Pell's
equation $v^2-du^2=4$.  Denote by $\Omega$ the set of all the positive
integers $d$ such that $d\equiv 0$ or 1 (mod 4) and such that $d$ is not
a square of an integer.

\proclaim{Theorem 1}    Let $N$ be a square free positive
integer, and let $n$ be a positive integer with $(n, N)=1$.
Put
$$L_n(s)=\sum_{m|N}\sum_{k|N}k^{1-2s}{\mu((m, k))\over (m,k)}
\sum_{d\in \Omega}\sum_u
\left({d\over m}\right){h_d\ln \epsilon_d\over (du^2)^s} $$
for $\Re s>1$, where the summation on $u$ is taken over
all the positive integers $u$ such that
$\sqrt{4n+d k^2 u^2}\in\Bbb Z$.   Then $L_n(s)$ is analytic for
$\Re s>1$ and can be extended by analytic continuation to the
half-plane $\Re s> 0$ except for a possible pole at $s=1/2$ and for
possible simple poles at $s=1, \frac 12\pm i\kappa_j$,
$j=1,2,\cdots$.  For any eigenvalue $\lambda>0$ of the Laplace-Beltrami
operator for $\Gamma_0(N)$, we have
$$\text{tr} T_n=2n^{i\kappa}\text{Res}_{s=1/2+i\kappa} L_n(s),$$
where $\kappa=\sqrt{\lambda-1/4}$.   \endproclaim

  The paper is organized as follows.  In section 2, we recall
the theory of Selberg's trace formula \cite{11}.  Elements
of the set $\Gamma^*$, which is defined in section 2, can be
divided into four types, the identity, hyperbolic, elliptic
and parabolic elements.  Next, in section 3 we
compute contributions of the identity, elliptic, hyperbolic
and parabolic elements to the trace formula.   A technical
part of this section (Lemma 3.3 -- Theorem 3.9) is to compute
the contribution from hyperbolic elements whose fixed points are
cusps of $\Gamma_0(N)$, and the result is given in Theorem 3.9.
By using the Selberg trace formula and by considering the
contributions to the trace formula
of the identity, elliptic , parabolic elements
and hyperbolic elements whose fixed points are cusps, we obtain
the analyticity information about a series formed by contributions
to the trace formula of hyperbolic elements whose fixed points
are not cusps, and a precise statement is given in Theorem 3.12.
In section 4, we compute explicitly the total contribution to
the trace formula of hyperbolic elements whose fixed points
are not cusps.  The result is stated in Theorem 4.6.
A technical part of this section is to relate
the number of certain indefinite primitive quadratic
forms, which are not equivalent under $\Gamma_0(N)$, to the
class number of indefinite primitive quadratic forms for
the full modular group $SL_2(\Bbb Z)$, and the relation is given
in Lemma 4.5.  Finally, the main theorem of this paper
follows from Theorem 3.12 and Theorem 4.6.

    The authors wish to thank Atle Selberg for his valuable
suggestions during preparation of the manuscript, and wish
to thank William Duke for the reference \cite{2}.

\heading
2.   The Selberg trace formula
\endheading

 Let $s$ be a complex number with $\Re s>1$.  Define
$$k(t)=(1+\frac t4)^{-s}$$
and
$$k(z, z^\prime)=k\left(\frac{|z-z^\prime|^2}{yy^\prime}\right),$$
for $z=x+iy$ and $z^\prime=x^\prime+iy^\prime$ in the upper half-plane.
Then $k(mz, mz^\prime)=k(z, z^\prime)$ for every $2\times 2$ matrix $m$ of
determinant one with  real entries.  The kernel $k(z, z^\prime)$
is of (a)-(b) type in the sense of Selberg \cite {11}, p.60.  Let
  $$ g(u)=\int_w^\infty k(t)\frac{dt}{\sqrt {t-w}}$$
with $w=e^u+e^{-u}-2$.   Write
  $$h(r)=\int_{-\infty}^{\infty} g(u) e^{iru}du. $$
Then
$$g(u)=\sqrt w\int_0^1(t+\frac w4)^{-s}t^{s-\frac 32}\frac{dt}
{\sqrt{1-t}}=c(1+\frac w4)^{\frac 12-s} \tag 2.1$$
 where $c=2\sqrt \pi\Gamma(s-\frac 12)\Gamma^{-1}(s)$.   Since
 $$\aligned h(r)&=c4^{s-\frac 12}\int_0^\infty
 (u+\frac 1u+2)^{\frac 12-s}u^{ir-1}du\\
&=c4^{s-{1\over 2}}\int_1^\infty (u+{1\over u}+2)^{{1\over 2}-s}
\left(u^{ir}+u^{-ir}\right){du\over u}\\
&=c{4^s(s-{1\over 2})\over (s-{1\over 2})^2+r^2}+A(r, s),\endaligned $$
where $A(r,s)$ is finite for $|\Im r|\leqslant 1/2$ and for
$\Re s>0$, we obtain that
$$\lim_{s\to 1/2+i\kappa} (s-{1\over 2}-i\kappa)h(r)=
\cases 4^{1/2+i\kappa}\sqrt \pi\frac {\Gamma(i\kappa)}{\Gamma(1/2+i\kappa)},
&\text{for $r=\pm\kappa$}; \\0, &\text{for $r\neq \pm\kappa$.}
\endcases\tag 2.2$$
{\it Remark.}  Atle Selberg told the second author that he
used the function $k(t)=(1+{t\over 4})^{-s}$ in some of his
unpublished works, which was convenient for computations
(cf. Selberg \cite {13}).

   Let $n$ be a positive integer with $(n, N)=1$.  Define
$$\Gamma^*=\cup_{\underset{0\leq b<d}\to{ad=n}}
\frac 1{\sqrt n}\left(\matrix d&{-b}\\0&a\endmatrix\right)\Gamma_0(N).$$
 Since $(n, N)=1$, we have that $T^{-1}\in\Gamma^*$ whenever
$T\in\Gamma^*$.    Every element of $\Gamma^*$ is represented
 uniquely in the form
$$\frac 1{\sqrt n}\left(\matrix d&{-b}\\0&a\endmatrix\right)\gamma$$
with $ad=n$, $0\leq b<d$ and $\gamma\in\Gamma_0(N)$.  It follows that
$\Gamma^*$ satisfies all the requirements given in \cite {11}, p.69.
 Let $\frak a_1, \frak a_2, \cdots, \frak a_{\nu(N)}$ with
$\nu(N)=\sum_{w|N, w>0}\varphi\left((w, {N\over w})\right)$
be a complete set of inequivalent cusps of $\Gamma_0(N)$,
where $\varphi$ is Euler's function.  We choose an element $\sigma_{\frak a_i}
\in PSL(2, \Bbb R)$  such that $\sigma_{\frak a_i}\infty=\frak a_i$ and
$\sigma_{\frak a_i}^{-1}\Gamma_{\frak a_i}
\sigma_{\frak a_i}=\Gamma_\infty$ for $i=1,2, \cdots, \nu(N)$.
The Eisenstein series $E_i(z, s)$ for the cusp $\frak a_i$ is defined by
    $$E_i(z, s)=\sum_{\gamma\in\Gamma_{\frak a_i}\backslash\Gamma_0(N)}
\left(\Im (\sigma_{\frak a_i}^{-1}\gamma z)\right)^s$$
for $\Re s>1$ when $z$ is in the upper half-plane.  Define
$$K(z, z^\prime)=\sum_{T\in \Gamma^*}k(z, T z^\prime ) $$
and
$$H(z, z^\prime)= \sum_{i=1}^{\nu(N)}\sum_{ad=n, \,0\leq b< d}
\frac 1{4\pi}\int_{-\infty}^\infty h(r)E_i(\frac{az+b}d, \frac 12+ir)
E_i(z^\prime, \frac 12 -ir) dr.$$
Denote by tr$_jT_n$ the trace of the Hecke operator $T_n$
acting on the space $E_{\lambda_j}$.  It follows
from (2.14) of \cite {11}, the argument of \cite {7}, pp.96-98,
Theorem 5.3.3 of
\cite {7}, and the spectral decomposition formula (5.3.12) of \cite {7} that
$$d(n)h(-\frac i2)+\sqrt n \sum_{j=1}^\infty h(\kappa_j)
\text{tr}_jT_n=\int_D \{K(z, z)-H(z, z)\}dz \tag 2.3$$
for $\Re s>1$, where $d(n)$ is the sum of positive divisors of $n$.

\heading
3.   Evaluation of components of the trace
\endheading

 For every element $T$ of $\Gamma^*$,
denote by $\Gamma_T$ the set of all the elements of $\Gamma_0(N)$ which commute
with $T$.  Put $D_T=\Gamma_T\backslash \Cal H$.
The elements of $\Gamma^*$ can be divided into four types,
of which the first consists of the identity element, while the
others are respectively the hyperbolic, the elliptic and the parabolic
elements.   If $T$ is not a parabolic element, put
$$c(T)=\int_{D_T} k(z, Tz)dz.$$

\subheading{ 3.1.   The identity component}

  If $\Gamma^*$ contains the identity element $I$, then
  $$c(I)=\int_{\Gamma_0(N)\backslash\Cal H}dz.$$

\subheading{ 3.2.   Elliptic components}

     There are only a finite number of elliptic conjugacy classes.

\proclaim{Lemma 3.1} Let $R$ be an elliptic element of $\Gamma^*$.  Then
$$c(R)=\frac \pi{2m\sin \theta}\int_0^\infty
\frac{k(t)}{\sqrt{t+4\sin^2\theta}}dt, $$
where $m$ is the order of a primitive element of $\Gamma_R$ and
where $\theta$ is defined by the formula trace$(R)=2\cos\theta$.
\endproclaim

  \demo{Proof}   Since $R$ is an elliptic element of $\Gamma^*$,
an element $\sigma\in PSL(2, \Bbb R)$ exists such that
$$\sigma R\sigma^{-1}=\left(\matrix {\cos\theta}&{-\sin\theta}\\
{\sin\theta}&{\cos\theta}\endmatrix\right)=\widetilde R$$
for some real number $0<\theta<\pi$.
Denote by $(\sigma\Gamma_0(N)\sigma^{-1})_{\widetilde R}$
the set of all the elements of $\sigma\Gamma_0(N)\sigma^{-1}$ which commute
with $\widetilde R$.  We have
$$c(R)=\int_{D_{\widetilde R}}k(z, \widetilde R z)dz$$
where $D_{\widetilde R}
=(\sigma\Gamma_0(N)\sigma^{-1})_{\widetilde R}\backslash \Cal H$.

  Let $\gamma=\left(\smallmatrix {\alpha}&{\beta}\\{\gamma}&{\delta}
\endsmallmatrix\right)$ be an element of $\Gamma_0(N)$ which has the same
fixed points as $R=\left(\smallmatrix a&b\\c&d\endsmallmatrix\right)$.
Then $(\alpha-\delta)c=\gamma(a-d)$ and $\beta c=\gamma b$.  It follows
that $\gamma$ commutes with $R$.  By Proposition 1.16 of \cite {14},
a primitive elliptic element $\gamma_0$ of $\Gamma_0(N)$ exists such that
$(\eta\Gamma\eta^{-1})_{\widetilde R}$ is generated by $\eta\gamma_0\eta^{-1}$.
Since $\eta\gamma_0\eta^{-1}$ commutes with $\widetilde R$, it is of the form
$$\left(\matrix {\cos\theta_0}&{-\sin\theta_0}\\
{\sin\theta_0}&{\cos\theta_0}\endmatrix\right)$$
for some real number $\theta_0$.   By Proposition 1.16 of \cite {14},
$\theta_0=\pi/m$ for some positive integer $m$.
 It follows from the argument of \cite {7}, p.99 that
$$c(R)=\frac 1m\int_0^\infty\int_{-\infty}^\infty
k\left(\frac{|z^2+1|^2}{y^2}\sin ^2\theta\right)dz.$$
By the argument of \cite {7}, p.100 we have
$$c(R)=\frac \pi{2m\sin \theta}\int_0^\infty
\frac{k(t)}{\sqrt{t+4\sin^2\theta}}dt. \qed$$
\enddemo

\subheading{ 3.3.   Hyperbolic components}

   Let $P$ be a hyperbolic element of $\Gamma^*$.  Then an element $\rho$
exists in $SL_2(\Bbb R)$ such that
$$\rho P\rho^{-1}=\left(\matrix {\lambda_P} &0\\
0&{\lambda_P^{-1}}\endmatrix\right)=\widetilde P$$
with $\lambda_P>1$.  The number $\lambda_P^2$ is called the norm of
$P$, and is denoted by $NP$.  It follows that
$$c(P)=\int_{D_{\widetilde P}}k(z, NP z)dz$$
where $D_{\widetilde P}
=(\rho\Gamma_0(N)\rho^{-1})_{\widetilde P}\backslash \Cal H$.
Let $P_0$ be a primitive hyperbolic element of $SL_2(\Bbb Z)$,
which generates the group of all elements of $SL_2(\Bbb Z)$
commutating with $P$.  Then there exists a hyperbolic element
$P_1\in\Gamma_0(N)$, which generates $\Gamma_P$, such that
$P_1$ is the smallest positive integer power of $P_0$
among all the generators of $\Gamma_P$ in $\Gamma_0(N)$.
Detail discussion about the ``primitive" hyperbolic
element $P_1$ is given in the proof of Theorem 4.6.

\proclaim{Theorem 3.2}  Let $P$ be a hyperbolic element of $\Gamma^*$
such that $\Gamma_P\neq \{1_2\}$.  If $P_1$ is a ``primitive'' hyperbolic
element of $\Gamma_0(N)$ which generates the group $\Gamma_P$, then
$$c(P)=\frac{\ln NP_1}{(NP)^{1/2}-(NP)^{-1/2}}g(\ln NP). $$
\endproclaim

\demo{Proof}   An argument similar to that made for the elliptic
elements shows that every element of $\Gamma_0(N)$, which has the same
fixed points as $P$, commutes with $P$.  Because $\rho P_1\rho^{-1}$
commutes with $\widetilde P$, it is of the form
$$\left(\matrix {\lambda_{P_1}}&0\\0&{\lambda_{P_1}^{-1}}\endmatrix\right)$$
for some real number $\lambda_{P_1}>1$.  Then
$$c(P)=\int_1^{NP_1}\frac{dy}{y^2}\int_{-\infty}^\infty
k\left(\frac{(NP-1)^2}{NP}\frac{|z|^2}{y^2}\right)dx.$$
The stated identity follows.  \qed\enddemo

  We next consider $c(P)$ for hyperbolic element $P\in\Gamma^*$ with
$\Gamma_P=\{1_2\}$.

\proclaim{Lemma 3.3}  Let $P={1\over \sqrt n}\left(\smallmatrix
A&B\\C&D\endsmallmatrix\right)$ be a hyperbolic element of $\Gamma^*$
such that $\Gamma_P=\{1_2\}$.  Then fixed points of $P$ are cusps
of $\Gamma_0(N)$.  Moreover, $\Gamma_P=\{1_2\}$
if and only if $A+D={1\over 2}(m+{4n\over m})$ for some divisor $m$ of
$4n$ with $m\neq 2\sqrt n$ for $C\neq 0$, and if and only if $A\neq D$
for $C=0$.
\endproclaim

 \demo{Proof}  Let $P$ a hyperbolic element of $\Gamma^*$
such that $\Gamma_P=\{1_2\}$.   If the fixed points of $P$ are not
rational numbers or infinity, then they are zeros of an irreducible polynomial
$ax^2+bx+c$ with $a\equiv 0$ mod$(N)$.  If $d=b^2-4ac$ and $(u, v)$ is
a solution of Pell's equation $v^2-du^2=4$, then
$$\left(\matrix {v-bu\over 2}&{-cu}\\{au}&{v+bu\over 2}\endmatrix\right)$$
belongs to $\Gamma_0(N)$ and has the same fixed points as $P$,
and hence it commutes with $P$.  This contradicts to $\Gamma_P=\{1_2\}$.
Therefore, fixed points of $P$ are rational numbers or infinity.
If $u/w$ is a fixed point of $P$, then
$$\left(\matrix {1+Nuw}&{-Nu^2}\\{Nw^2}&{1-Nuw}\endmatrix\right)$$
is a parabolic element of $\Gamma_0(N)$ and has $u/w$ as
a fixed point.  Hence, fixed points of $P$ are cusps of $\Gamma_0(N)$.

  Conversely, if $T$ is an element of $\Gamma^*$ having two distinct
fixed points with at least one of them being a rational number,
then $\Gamma_T=\{1_2\}$.   Otherwise, let
$r_1, r_2$ be the fixed points of $T$, and let $\gamma\in\Gamma_0(N)$
($\gamma\neq 1_2$) with $\gamma T=T\gamma$.  Then either $\gamma$ or
$\gamma^2$ ($\neq 1_2$), say $\gamma^2$, has $r_1, r_2$ as its fixed points,
and hence it is a hyperbolic element of $\Gamma_0(N)$.  Since
rational points are cusps of $\Gamma_0(N)$, by Proposition 1.17 of
\cite {14} $\gamma^2$ is parabolic.  This is a contradiction, and hence
we have $\Gamma_T=\{1_2\}$.  Thus, we have proved that $\Gamma_T=\{1_2\}$
for an element $T={1\over \sqrt n}
\left(\smallmatrix a&b\\c&d\endsmallmatrix\right)\in\Gamma^*$
if and only if $a+d={1\over 2}(m+{4n\over m})$ for some divisor $m$ of
$4n$ with $m\neq 2\sqrt n$ if $c\neq 0$ or $a\neq d$ if $c=0$.
\qed\enddemo

 Every cusp of $\Gamma_0(N)$ is equivalent to one of
the following inequivalent cusps
$$ {u\over w}\,\,\text{with}\,\, u, w>0,
\,\, (u, w)=1,\,\, w|N.\,\,\tag 3.1$$
Two such cusps $u/w$ and $u_1/w_1$ are $\Gamma_0(N)$-equivalent
if and only if $w=w_1$ and $u\equiv u_1$ modulo $(w, N/w)$.
Let $\frak a=u/w$ be given as in (3.1).  By (2.2) and (2.3)
of Deshouillers and Iwaniec \cite {1}, we have
$$\Gamma_{\frak a}=\left\{\left(\matrix {1+cu/w}&{-cu^2/w^2}\\c&{1-cu/w}
\endmatrix\right):\,\, c\equiv 0
\,\,(\text{mod}\,[w^2, N])\right\} \tag 3.2$$
and
$$\sigma_{\frak a}\infty=\frak a\,\,\,\,\text{and}\,\,\,\,
\sigma_{\frak a}^{-1}\Gamma_{\frak a}\sigma_{\frak a}=\Gamma_\infty$$
where
$$\sigma_{\frak a}=\left(\matrix {\frak a\sqrt{[w^2, N]}}&0\\
{\sqrt{[w^2, N]}}&{1/\frak a \sqrt{[w^2, N]}}\endmatrix\right).
\tag 3.3$$
Let $P$ be a hyperbolic element of $\Gamma^*$
such that $\Gamma_P=\{1_2\}$.  Assume that
$\frak a$ is a fixed point of $P$.  Then $\infty$ is a fixed point
of $\sigma_{\frak a}^{-1}P\sigma_{\frak a}$, and hence there exist
positive numbers $a,d$ with $ad=n, a\neq d$ such that $P$ is of the form
$$P={1\over\sqrt n}\sigma_{\frak a}\left(\matrix a&b\\0&d\endmatrix\right)
\sigma_{\frak a}^{-1}={1\over \sqrt n}\left(\matrix{a-b[w^2, N]
{u\over w}}&{b[w^2, N]{u^2\over w^2}}\\
{(a-d){w\over u}-b[w^2, N]}&{d+b[w^2, N]{u\over w}}
\endmatrix\right).\tag 3.4$$
We claim that $a, d$ are integers.  Since $P={1\over\sqrt n}\left(\smallmatrix
A&B\\C&D\endsmallmatrix\right)\in\Gamma^*$, we have $A,B,C, D\in\Bbb Z$,
$AD-BC=n$, and $C\equiv 0$ mod($N$).  By (3.4), we find that
$$a=D+C{u\over w},\,\, d=D-B{w\over u},\,\, b=B{w^2\over u^2[w^2, N]}\tag 3.5$$
and
$$A-D=C{u\over w}-B{w\over u}.\tag 3.6$$
Since $(w, u)=1$, we have $u|B$ by (3.6).  It follows from (3.5)
that $a, d$ are integers.  Since $c(P)$ depends only on the conjugacy class
$\{P\}$ represented by $P$, we can replace $P$
 by $\gamma^{-1}P\gamma$ without changing the value of $c(P)$.
Replacing $P$ by $\gamma^{-1}P\gamma$ for some element $\gamma\in
\Gamma_{\frak a}$, we can assume without loss of generality that
$0\leqslant b<|a-d|$ in (3.4).

  \proclaim{Lemma 3.4}  Let $\frak a=u/w$ be given as in (3.1).  Then
$$P={1\over\sqrt n}\sigma_{\frak a}\left(\matrix a&b\\0&d\endmatrix\right)
\sigma_{\frak a}^{-1} $$
with $a, d\in\Bbb Z, ad=n, a\neq d$ is a hyperbolic element of $\Gamma^*$ with
$\Gamma_P=\{1_2\}$ and $P(\frak a)=\frak a$ if and only if
$(w, N/w)|(a-d)$ with $b$ being chosen so that
$(a-d){w\over u}-b[w^2, N]$ is divisible by $N$.\endproclaim

  \demo{Proof}  Let $P$ be a hyperbolic element of $\Gamma^*$ with
$\Gamma_P=\{1_2\}$ and $P(\frak a)=\frak a$.
By (3.4), $(a-d){w\over u}-b[w^2, N]$ is divisible
by $N$, and there exists an integer $k$ such that ${u\over w}N|a-d-kw$.
This implies that $(w, N/w)|a-d$, which
can also be seen directly from the identities $P(\frak a)=\frak a$
and $a+d={1\over 2}(m+{4n\over m})$ for some divisor $m$ of $4n$.

   Conversely, let $a, d\in\Bbb Z^+, ad=n, a\neq d$ and
$(w, N/w)|a-d$.  Write $a-d=\ell (w, N/w)$ for some integer $\ell$.
Since $uN/w(w, N/w)$ and $w/(w, N/w)$ are coprime, there exist integers
$\lambda$ and $\tau$ such that $\lambda\frak a N+\tau w=(w, N/w)$,
and hence $\ell\lambda \frak a N=a-d-\ell\tau w$.  This means that
there exists an integer $k$ such that ${u\over w}N|a-d-kw$.  If
we choose $A=a-kw, B=ku, C={w\over u}(a-d-kw)$ and $D=d+kw$,
then $AD-BC=n$ and $C\equiv 0$ mod($N$).  Let
$P={1\over\sqrt n}\left(\smallmatrix A&B\\C&D\endsmallmatrix\right)$.
Then $P$ can be expressed as in (3.4), $P(\frak a)=\frak a$, and
$\Gamma_P=\{1_2\}$.  Let $\gamma=C/(C, D)$ and
$\delta=D/(C, D)$.  There exist integers $\alpha, \beta$ such that
$\alpha\delta-\beta\gamma=1$.  Since $(n, N)=1$, we have
$\gamma\equiv 0$ mod($N$).  It follows that
$$P={1\over \sqrt n}\left(\matrix{n\over C,D)}&{B\alpha-A\beta}\\0&{(C,D)}\endmatrix\right)\left(\matrix
\alpha &\beta\\ \gamma &\delta\endmatrix\right).$$
This implies that $P\in\Gamma^*$.  Thus, we have proved that
$$P={1\over\sqrt n}\sigma_{\frak a}\left(\matrix a&b\\0&d\endmatrix\right)
\sigma_{\frak a}^{-1} $$
with $a, d\in\Bbb Z, ad=n, a\neq d$ is an element of $\Gamma^*$ with
$\Gamma_P=\{1_2\}$ and $P(u/w)=u/w$ if and only if there exists an integer
$k$ such that ${u\over w}N|a-d-kw$, that is, if and only if
$(w, N/w)|(a-d)$ with $b$ being chosen so that
$(a-d){w\over u}-b[w^2, N]$ is divisible by $N$.
\qed\enddemo

  \proclaim{Lemma 3.5}   Let $\frak a=u/w$ be given as in (3.1).
Let $a, d\in\Bbb Z^+, ad=n, a\neq d$.  If $(w, N/w)|a-d$, then there
are exactly $|a-d|$ number of $\Gamma_0(N)$-inequivalent hyperbolic
elements $P\in\Gamma^*$ with $P(\frak a)=\frak a$ and $\Gamma_P
=\{1_2\}$, which are of the form
$$P={1\over\sqrt n}\sigma_{\frak a}\left(\matrix a&b\\0&d\endmatrix\right)
\sigma_{\frak a}^{-1}$$
with $0\leqslant b<|a-d|$.  \endproclaim

  \demo{Proof}    For $a, d\in\Bbb Z, ad=n, a\neq d$, let
$$P={1\over\sqrt n}\sigma_{\frak a}\left(\matrix a&b\\0&d\endmatrix\right)
\sigma_{\frak a}^{-1}\,\,\,\,\text{and}\,\,\,\,
P^\prime={1\over\sqrt n}\sigma_{\frak a}\left(\matrix a&b^\prime\\0&d
\endmatrix\right) \sigma_{\frak a}^{-1}$$
be two elements of $\Gamma^*$ with $0\leqslant b, b^\prime<|a-d|$
such that $\Gamma_P=\{1_2\}, P(\frak a)=\frak a$
and $\Gamma_{P^\prime}=\{1_2\}, P^\prime(\frak a)=\frak a$.
We claim that $b^\prime=b+\ell$ for some integer $\ell$.
In fact, for $P$ there exists an integer $k$ such that
${u\over w}N|a-d-kw$, and for $P^\prime$ there exists an integer
$k^\prime$ such that ${u\over w}N|a-d-k^\prime w$.  Then we have
${u\over w}N|(k^\prime-k)w$, that is,
$$k^\prime=k+{\ell uN\over w(w, N/w)}$$
for some integer $\ell$.  By using (3.4), we find that
$$b^\prime=b+{\ell wN\over [w^2, N](w, N/w)}=b+\ell.$$
The claim then follows.  Conversely, given an integer $\ell$ with
$0\leqslant b+\ell< |a-d|$, let $b^\prime=b+\ell$ and
$$k^\prime=k+{\ell uN\over w(w, N/w)}.$$
Put $A^\prime=a-k^\prime w, B^\prime=k^\prime u,
C^\prime={w\over u}(a-d-k^\prime w)$ and
$D^\prime=d+k^\prime w$.  Then
$$P^\prime={1\over\sqrt n}\sigma_{\frak a}\left(\matrix
a&b^\prime\\0&d\endmatrix\right) \sigma_{\frak a}^{-1}$$
is a hyperbolic element of $\Gamma^*$,
$\Gamma_{P^\prime}=\{1_2\}$, and $P^\prime(\frak a)=\frak a$.
Moreover, if $b^\prime\neq b$ modulo $|a-d|$, then
$P^\prime$ and $P$ are not $\Gamma_0(N)$-equivalent.
Otherwise, an element $\gamma\in\Gamma_0(N)$ exists such that
$P^\prime\gamma=\gamma P$.   Put $\gamma=\left(\smallmatrix
*&*\\C&*\endsmallmatrix\right)$ and $\sigma_{\frak a}^{-1}
\gamma\sigma_{\frak a}=\left(\smallmatrix\alpha &\beta\\
\eta &\delta\endsmallmatrix\right)$.  By using the relation
$P^\prime\gamma=\gamma P$, we obtain that
$(a-d)\beta=\alpha b-\delta b^\prime$ and $\eta=0$.
Since $\gamma\in\Gamma_0(N)$, by using $\eta=0$ and
the relation $\sigma_{\frak a}^{-1}
\gamma\sigma_{\frak a}=\left(\smallmatrix\alpha &\beta\\
0&\delta\endsmallmatrix\right)$ we find that $\alpha$
and $\delta$ are integers, and hence $\alpha=\delta=1$.
Furthermore, we obtain that
$\beta=-C/[w^2, N]$.  Since $\sigma_{\frak a}^{-1}
\gamma\sigma_{\frak a}=\left(\smallmatrix 1&\beta\\
0&1\endsmallmatrix\right)$, we have $(\sigma_{\frak a}^{-1}\gamma
\sigma_{\frak a})(\infty)=\infty$, that is,
$\gamma (\frak a)=\frak a$.  Hence $\gamma\in\Gamma_{\frak a}$,
which is given by (3.2).  This implies that $\beta$ is
an integer.  It follows that $(a-d)|(b-b^\prime)$.
This is a contradiction.
Therefore, for fixed $a, d\in\Bbb Z^+, ad=n, a\neq d$, there
are exactly $|a-d|$ $\Gamma_0(N)$-inequivalent hyperbolic
elements $P$ in $\Gamma^*$ with $P(u/w)=u/w$ and $\Gamma_P
=\{1_2\}$, which are of the form
$$P={1\over\sqrt n}\sigma_{\frak a}\left(\matrix a&b\\0&d\endmatrix\right)
\sigma_{\frak a}^{-1}$$
with $0\leqslant b<|a-d|$.  \qed\enddemo

  \proclaim{Lemma 3.6}  Let $\frak a=u/w$, and let
$$P={1\over\sqrt n}\sigma_{\frak a}\left(\matrix a&b\\0&d\endmatrix\right)
\sigma_{\frak a}^{-1}\,\,\,\,\text{and}\,\,\,\,
P^\prime={1\over\sqrt n}\sigma_{\frak a}\left(\matrix d&b^\prime\\0&a
\endmatrix\right) \sigma_{\frak a}^{-1}$$
be two hyperbolic elements of $\Gamma^*$ with  $P(\frak a)=\frak a,
P^\prime(\frak a)=\frak a$, $\Gamma_P=\{1_2\}$ and $\Gamma_{P^\prime}
=\{1_2\}$.   Then $P$ is $\Gamma_0(N)$-conjugate to $P^\prime$ for some
number $b^\prime$ if and only if the two fixed points of $P$
are $\Gamma_0(N)$-conjugate.   That is, $P$ is $\Gamma_0(N)$-conjugate
to $P^\prime$ for some number $b^\prime$ if and only if
$$\left({(d-a){w\over u}+b[w^2, N]\over \ell}, N\right)=w,$$
where $\ell$ is the greatest common divisor of $b[w^2, N]{u\over w}$
and $(d-a){w\over u}+b[w^2, N]$. \endproclaim

 \demo{Proof}   We first assume that the two fixed points of $P$
are $\Gamma_0(N)$-conjugate.  That is, there exists an
element $T=\left(\smallmatrix \alpha&\beta\\
\gamma&\delta\endsmallmatrix \right)\in\Gamma_0(N)$, which is
not the identity, such that $\frak a, T(\frak a)$ are the two
fixed points of $P$.   Let $P^\prime=T^{-1}PT$.  Then
$P^\prime(\frak a)=\frak a$.  Since $\Gamma_P=\{1_2\}$,
we have $\Gamma_{P^\prime}=\{1_2\}$.  We claim that
$P^\prime$ is of the form given in the statement of the lemma.

  We can write $P^\prime$ in the form
$$P^\prime={1\over\sqrt n}\sigma_{\frak a}\left(\matrix
{a^\prime}&{b^\prime}\\0&{d^\prime}\endmatrix\right)
\sigma_{\frak a}^{-1}$$
with $a^\prime, d^\prime\in\Bbb Z^+$ and
$a^\prime d^\prime=n$.  Replacing $T$ by $T\gamma$
for some element $\gamma\in \Gamma_{\frak a}$,
we assume that $0\leqslant b^\prime<|a^\prime-d^\prime|$.
Since $P^\prime=T^{-1}PT$, we have
$$(\sigma_{\frak a}^{-1}T\sigma_{\frak a})\left(\matrix
a^\prime&b^\prime\\0&d^\prime\endmatrix\right)
=\left(\matrix a&b\\0&d\endmatrix\right)
(\sigma_{\frak a}^{-1}T\sigma_{\frak a}). $$
In general, from the identity
$$\left(\matrix A&B\\C&D\endmatrix\right)
\left(\matrix {a^\prime}&{b^\prime}\\0&{d^\prime}\endmatrix\right)
=\left(\matrix a&b\\0&d\endmatrix\right)
\left(\matrix A&B\\C&D\endmatrix\right)$$
with $AD-BC\neq 0$, we obtain that
$(d-a^\prime)C=0, (d-d^\prime)D=b^\prime C$ and
$(a^\prime-a)A=bC$.
If $C=0$, we have $a^\prime=a$ and $d^\prime=d$, and
hence $P^\prime$ is of the form
$$P^\prime={1\over\sqrt n}\sigma_{\frak a}\left(\matrix
a&{b^\prime}\\0&d\endmatrix\right)\sigma_{\frak a}^{-1}.$$
Since $P^\prime$ is $\Gamma_0(N)$-conjugate to $P$ and
$0\leqslant b^\prime<|a-d|$, we must have
$P^\prime=P$ by the proof of Lemma 3.5.  That is,
$T\in\Gamma_P$, and hence $T=1_2$.
 A contradiction is then derived.  Therefore, we
must have $C\neq 0$.  It follows that $a^\prime=d$.
Since $ad=n$ and $a^\prime d^\prime=n$, we have
$d^\prime=a$, and our claim follows.

  Conversely, suppose that there exists an element
 $T=\left(\smallmatrix A&B\\C&D\endsmallmatrix\right)\in
\Gamma_0(N)$ such that $T^{-1}PT=P^\prime$, that is,
$$\sigma_{\frak a}^{-1}T^{-1}\sigma_{\frak a}
\left(\matrix a&b\\0&d\endmatrix\right)
=\left(\matrix d&{b^\prime}\\0&a\endmatrix\right)
\sigma_{\frak a}^{-1}T^{-1}\sigma_{\frak a}.\tag 3.7$$
If we write
$$\sigma_{\frak a}^{-1}T^{-1}\sigma_{\frak a}=\left(
\matrix{D-B{w\over u}}&{-B w^2\over [w^2, N]u^2}\\
{{u\over w}[w^2, N](A+B{w\over u}-C{u\over w}-D)}
&{A+B{w\over u}}\endmatrix\right),$$
then (3.7) is equivalent to the system of equations
 $$\cases (a-d)(A+B{w\over u})
=b[w^2, N]{u\over w}(A+B{w\over u}
-C{u\over w}-D),\\
 b(D-B{w\over u})=b^\prime(A+B{w\over u}).
\endcases $$
Write this system of equations as
$$\cases\left((a-d){w\over u}-b[w^2, N]\right)u(A+B
{w\over u})=b[w^2, N]{u\over w}(-C u-D w),\\
b(D-B{w\over u})=b^\prime(A+B{w\over u}).
\endcases \tag 3.8$$
The first equation of (3.8) can be written as
$$T(\frak a)={\frak a\over 1-{a-d\over b[w^2, N]\frak a}},$$
which is the second fixed point of $P$.  That is,
the two fixed points of $P$ are $\Gamma_0(N)$-conjugate.
Since $\frak a={u\over w}$ and
$${b[w^2, N]u/w\over (d-a){w\over u}+b[w^2, N]}$$
are the two fixed points of $P$,
the last statement of the lemma follows from Lemma 3.6
of Deshouillers and Iwaniec \cite {1}. \qed\enddemo

\proclaim{Lemma 3.7}  Let $a,d\in \Bbb Z^+, a\neq d, ad=n$,
and let $\frak a=u/w$ be given as in (3.1) with
$(w, N/w)|(a-d)$.  Assume that $\frak a$ and
$\frak a^\prime=u_1/w_1$ are the two
distinct fixed point of a hyperbolic element $P\in\Gamma^*$.
Then $(w^\prime, N/w^\prime)|(a-d)$
for $w^\prime=(w_1, N)$.  In other words, if
$P={1\over\sqrt n}\sigma_{\frak a}\left(\smallmatrix a&b
\\0&d\endsmallmatrix\right)\sigma_{\frak a}^{-1}$
is a hyperbolic element of $\Gamma^*$ with
 $P(\frak a^\prime)=\frak a^\prime$, then it can also
be written in the form
$$P={1\over\sqrt n}\sigma_{\frak a^\prime}\left(\matrix d&b^\prime
\\0&a\endmatrix\right)\sigma_{\frak a^\prime}^{-1}$$
for some number $b^\prime$. \endproclaim

  \demo{Proof}  It follows from the argument made in the
paragraph preceding Lemma 3.4 that there exist positive
integers $a^\prime, d^\prime$ with $a^\prime d^\prime=n$ such that
$$P={1\over\sqrt n}\sigma_{\frak a^\prime}\left(\matrix a^\prime
&b^\prime \\0&d^\prime\endmatrix\right)\sigma_{\frak a^\prime}^{-1}$$
for some number $b^\prime$.  Then we have
$$\sigma_{\frak a}^{-1}\sigma_{\frak a^\prime}\left(\matrix a^\prime
&b^\prime \\0&d^\prime\endmatrix\right)
=\left(\matrix a&b\\0&d\endmatrix\right)
\sigma_{\frak a}^{-1}\sigma_{\frak a^\prime}. \tag 3.9$$
Write $\sigma_{\frak a}^{-1}\sigma_{\frak a^\prime}=\left(
\smallmatrix A&B\\C&D\endsmallmatrix\right)$.  It
follows from (3.9) that $(d-a^\prime)C=0, (a^\prime -a)A=bC$
and $(d-d^\prime)D=b^\prime C$.  Since $\frak a\neq\frak a^\prime$,
by using (3.3) we find that $C\neq 0$, and hence we have
$a^\prime=d$.  It follows that $d^\prime=a$, that is,
$P$ can also be written in the form stated in the lemma.
Then by the argument made in the first paragraph of the
proof of Lemma 3.4, we have $(w^\prime, N/w^\prime)|(a-d)$.
 \qed\enddemo

   For a large positive number $Y$, define
$$D_Y=\{z\in D: \Im \sigma_{\frak a_i}^{-1}z<Y,\,\, i=1,2,\cdots,h\}.$$
Let
$$(D_P)_Y=\bigcup_{\gamma\in\Gamma_0(N)} \gamma D_Y.$$
Write
$$c(P)_Y=\int_{(D_P)_Y}k(z, Pz)dz. $$

 \proclaim{Lemma 3.8}  Let $\frak a=u/w$ be given as in (3.1), and let
$P={1\over\sqrt n}\sigma_{\frak a}\left(\smallmatrix a&b
\\0&d\endsmallmatrix\right)\sigma_{\frak a}^{-1}$
with $a, d\in\Bbb Z^+, ad=n, a\neq d$ be a hyperbolic element
of $\Gamma^*$ with $\Gamma_P=\{1_2\}$ and $P(\frak a)=\frak a$.
Then
$$c(P)_Y=\frac{\sqrt n\ln \{(a-d)^2[w^2, N]Y/2\rho\}}{|a-d|}
g(\ln \frac ad) +\int_1^\infty  k\left(\frac{(a-d)^2}nt\right)
\frac{\ln t}{\sqrt {t-1}}dt+o(1) $$
where  $o(1)\to 0$ as $Y\to\infty$ and
$\rho=C^2/2[C^2/\ell^2, N]Y$ with $C=(a-d){w\over u}-b[w^2, N]$
and with $\ell$ being given as in Lemma 3.6.  \endproclaim

  \demo{Proof}   Let
$$\gamma=\left(\matrix p&r\\q&s\endmatrix \right)$$
be an element of $SL_2(\Bbb R)$.  The linear fractional
transformation, which takes
every complex $z$ in the upper half-plane into $\gamma (z)$, maps the
horizontal line $\Im z=Y$ into a circle of radius $\frac 1{2q^2 Y}$
with center at $\frac pq+\frac i{2q^2 Y}$.
Let $$\mu=\left(\matrix 1&{\frac b{a-d}}\\0&1\endmatrix \right).$$
Then
$$P={1\over\sqrt n}\sigma_{\frak a}\mu^{-1}\left(\matrix a&0\\0&d
\endmatrix\right)\mu\sigma_{\frak a}^{-1}.$$
Note that
$$\mu\sigma_{\frak a}^{-1}=\left(\matrix {{w/u\over\sqrt{[w^2, N]}}
-{b\sqrt{[w^2, N]}\over a-d}}&{bu\sqrt{[w^2, N]}\over (a-d)w}
\\{-\sqrt{[w^2, N]}}&{{u\over w}\sqrt{[w^2, N]}}\endmatrix\right).$$
Since
$$P={1\over\sqrt n}\sigma_{\frak a}\mu^{-1}\left(\matrix {1/y}&0\\
0&y\endmatrix\right)\left(\matrix a&0\\0&d \endmatrix\right)
\left(\matrix y&0\\0&{1/y}\endmatrix\right)\mu\sigma_{\frak a}^{-1}$$
for any $y\neq 0$, by choosing $y=(a-d)\sqrt{[w^2, N]}$ we obtain that
 $$c(P)_Y=\int_{\mu_{\frak a}\{(D_P)_Y\}} k(z, {a\over d}z)dz$$
where
$$\mu_{\frak a}=\left(\matrix {(a-d){w\over u}-b[w^2, N]}
&{b[w^2, N]{u\over w}}\\{-1\over a-d}&{u\over (a-d)w}
\endmatrix\right).$$

The linear transformation $z\to\sigma_{\frak a} z$ maps the
half-plane $\Im z>Y$ into a disk $D_{\frak a}$ of radius
${1\over 2[w^2, N]Y}$ with center at $\frak a+{i\over 2[w^2, N]Y}$.
Then the transformation $z\to \mu_{\frak a}(z)$
maps the disk $D_{\frak a}$ into the half-plane
$\Im z> (a-d)^2[w^2, N]Y$.
Let
$$\frak p={-b[w^2, N]u/w\over (a-d){w\over u}-b[w^2, N]}.$$
If
$$\sigma_{\frak p}=\left(\matrix {\frak p\sqrt
{[C^2/\ell^2, N]}}&0\\ {\sqrt{[C^2/\ell^2, N]}}&
{1/\frak p\sqrt{[C^2/\ell^2, N]}}\endmatrix\right),$$
 then $\sigma_{\frak p}\infty=\frak p$
and $\sigma_{\frak p}^{-1}\Gamma_{\frak p}\sigma_{\frak p}
=\Gamma_{\infty}$.  By the definition of $D_Y$,
 the image of the half-plane $\Im z>Y$ under the
linear transformation $z\to\sigma_{\frak p}z$ is not
contained in $(D_P)_Y$.   Since the linear transformation
$z\to (\mu_{\frak a}\sigma_{\frak p})(z)$ maps the
half-plane $\Im z>Y$ into a disk $D_\rho$ of radius
$\rho$ centered at $i\rho$, where $\rho=C^2/2[C^2/\ell^2, N]Y$,
the disk $D_\rho$ is not contained in
$\mu_{\frak a}\{(D_P)_Y\}$.  It follows that
$$\aligned c(P)_Y &=\int_0^\pi d\theta
\int_{2\rho \sin \theta}^{(a-d)^2[w^2, N]Y/\sin \theta}k
\left({(a-d)^2\over n\sin^2\theta}
\right)\frac{dr}{r\sin^2\theta}+o(1)\\
&=\int_1^\infty k\left(\frac{(a-d)^2}nt\right)
\frac{\ln \{(a-d)^2[w^2, N]Yt/2\rho\}}{\sqrt {t-1}}dt
+o(1) \endaligned $$
where $o(1)$ has a limit zero when $Y\to\infty$, and hence we have
$$c(P)_Y=\frac{\sqrt n\ln \{(a-d)^2[w^2, N]Y/2\rho\}}{|a-d|}
g(\ln \frac ad) +\int_1^\infty  k\left(\frac{(a-d)^2}nt\right)
\frac{\ln t}{\sqrt {t-1}}dt+o(1). $$
This completes the proof of the lemma. \qed\enddemo

  In the rest of the paper, we shall indicate explicitly when
we assume that $N$ is square free.

\proclaim{Theorem 3.9}  Let $N$ be a square free positive integer.
  Then we have
$$\aligned \sum_{\{P\}, \Gamma_P=\{1_2\}}& c(P)_Y
=\nu(N)\sqrt n\sum_{ad=n, d>0, a\neq d}g(\ln{a\over d})\ln Y\\
&+{\sqrt n\over 2}\sum_{w|N, w>0}\sum_{ad=n, d>0, a\neq d}
\sum_b{\ln \{(a-d)^2wN[{C^2\over\ell^2},N]/C^2\}
\over |a-d|}g(\ln \frac ad)\\
&+\sum_{ad=n, d>0,a\neq d}{1\over 2}\nu(N)
 |a-d|\int_1^\infty  k\left(\frac{(a-d)^2}nt\right)
\frac{\ln t}{\sqrt {t-1}}dt+o(1)\endaligned$$
with $C=(a-d)w-bwN$, where $\ell=(C, bN)$ and $\sum_b$
is over the $|a-d|$ numbers such that $C, bN/w\in \Bbb Z, N|C$
 and $0\leqslant b<|a-d|$. \endproclaim

  \demo{Proof}  Let $P$ and $P^\prime$ be hyperbolic elements
in $\Gamma^*$ such that $\Gamma_P=\{1_2\}$ and
$\Gamma_{P^\prime}=\{1_2\}$.  Assume that $P^\prime$ and
$P$ are $\Gamma_0(N)$-conjugate, that is, $T^{-1}P^\prime T=P$
for some element $T\in\Gamma_0(N)$.  If $\frak a_1, \frak a_2$
are the two fixed points of $P$, then $T(\frak a_1), T(\frak a_2)$
are the two fixed points of $P^\prime$.  In other words,
if at least one of the two fixed points of $P$ is not
$\Gamma_0(N)$-conjugate to a fixed point of $P^\prime$,
then $\{P^\prime\}$ and $\{P\}$ represent different
$\Gamma_0(N)$-conjugacy classes.

   Denote $P(a, d;b;w)={1\over\sqrt n}\sigma_{1/w}\left(
\smallmatrix a&b\\0&d\endsmallmatrix\right)\sigma_{1/w}^{-1}$
with $\sigma_{1/w}$ being given as in (3.3).
It follows from Lemma 3.4, Lemma 3.5, Lemma 3.6 and Lemma 3.7
that
$$\sum_{w|N, w>0}\sum_{ad=n, d>0, a\neq d}
\sum_b c\left(P(a,d;b;w)\right)_Y
=2\sum_{\{P\}, \Gamma_P=\{1_2\}} c(P)_Y.$$
The stated identity then follows from Lemma 3.8.  \qed\enddemo

\subheading{ 3.4.   Parabolic components}

  Let $S$ be a parabolic element of $\Gamma^*$.   An argument similar to
that made for the elliptic elements shows that
every element of $\Gamma_0(N)$ which has the same fixed point as
$S$ commutes with $S$.   If $\frak a=u/w$ is the fixed point
of $S$, then $\Gamma_S=\Gamma_{\frak a}$, and hence we have
$\sigma_{\frak a}^{-1}\Gamma_S \sigma_{\frak a}=\Gamma_\infty$,
where $\sigma_{\frak a}$ is given as in (3.3).
Since $\sigma_{\frak a}^{-1} S\sigma_{\frak a}$ commutes
with every element of $\Gamma_\infty$, we have
$$\sigma_{\frak a}^{-1} S\sigma_{\frak a}={1\over \sqrt n}
\left(\matrix a&b\\0&a\endmatrix\right)$$
for some real numbers $a, b$ with $a^2=n$.
An argument similar to that made in
the paragraph following (3.4) shows that $a$ is an integer.
This implies that  $\Gamma^*$ has  parabolic elements only
if $n$ is the square of an integer.
Furthermore,  by (3.4) we see that elements of the form
$$S=\sigma_{\frak a}\left(\matrix 1&{b/{\sqrt n}}\\0&1\endmatrix\right)
\sigma_{\frak a}^{-1}, \, \, \, \, \, \, \,0\neq b\in \Bbb Z$$
constitute a complete set of representatives for the conjugacy
classes of parabolic elements of $\Gamma^*$ having $\frak a$ as
its fixed point.  It follows that
$$\sum_{\{S\}}\int_{D_Y}k(z, Sz)dz =\int_0^Y\int_0^1
\sum_{0\neq b\in \Bbb Z}k(z, z+\frac b{\sqrt n})dz+o(1), $$
where the summation on $\{S\}$ is taken over all parabolic
classes represented by parabolic elements whose fixed point is
$\frak a=u/w$ and where $o(1)$ tends to zero as $Y\to\infty$.
Define $\delta_n$ to be one if $n$ is the square of an
integer and to be zero otherwise.

 \proclaim{Theorem 3.10}  Put
$$c(\infty)_Y=\delta_n\nu(N)\int_0^Y\int_0^1\sum_{0\neq b\in \Bbb Z}
k(z, z+\frac b{\sqrt n})dz -\int_{D_Y} H(z, z)dz.$$
Then
$$\aligned \frac {c(\infty)_Y}{\sqrt n} &=\delta_n\nu(N)
g(0)\ln{\sqrt n\over 2}
-\nu(N)\sum_{ad=n, d>0, a\neq d}g(\ln{a\over d})\ln Y\\
&-{\delta_n\nu(N)\over 2\pi}\int_{-\infty}^\infty
h(r)\frac{\Gamma^\prime}\Gamma(1+ir)dr +{1\over 4} h(0)
\{\delta_n\nu(N)+d(n)\sum_{i=1}^{\nu(N)}
\varphi_{ii}({1\over 2})\}\\
&+\sum_{i,j=1}^{\nu(N)}\sum_{ad=n, d>0}
{1\over 4\pi}\int_{-\infty}^\infty h(r)
\left({a\over d}\right)^{ir}\varphi_{ij}^\prime({1\over 2}+ir)
\varphi_{ij}({1\over 2}-ir)dr+o(1)\endaligned$$
where $o(1)$ tends to zero as $Y\to\infty$ and
$\varphi_{ii}(1/2)\neq \infty$ (cf. Selberg \cite {13}).
\endproclaim

 \demo{Proof}  By the argument of \cite {7}, pp.102--106 we have
$$\aligned &\frac 1{\sqrt n}\int_0^Y \int_0^1 \sum_{0\neq b\in \Bbb Z}
k(z, z+\frac b{\sqrt n})dz \\
& =g(0)\ln(\sqrt n Y)-\frac 1{2\pi}\int_{-\infty}^\infty
h(r)\frac{\Gamma^\prime}\Gamma(1+ir)dr -g(0)\ln 2
+\frac 14 h(0)+o(1).\endaligned $$
Let
$$\varphi_{ij, m}(s)=\sum_c{1\over |c|^{2s}}\left(\sum_d
e(md/c)\right)$$
where summations are taken over $c>0, d$ modulo $c$ with
$\left(\smallmatrix *&*\\ c&d\endsmallmatrix\right)\in
\sigma_{\frak a_i}^{-1}\Gamma_0(N)\sigma_{\frak a_j}$.
Then we have
$$\aligned E_i(\sigma_{\frak a_j}z, s)&=\delta_{ij}y^s+
\varphi_{ij}(s)y^{1-s}\\
&+{2\pi^s\sqrt y\over\Gamma(s)}\sum_{m\neq 0}
|m|^{s-{1\over 2}}K_{s-{1\over 2}}(2|m|\pi y)
\varphi_{ij, m}(s)e(mx)\endaligned$$
where
$$\varphi_{ij}(s)={\sqrt\pi\Gamma(s-{1\over 2})\over
\Gamma(s)}\varphi_{ij, 0}(s).$$
By (1.1) and the Maass-Selberg relation (Theorem 2.3.1 of \cite {7}),
we obtain that
$$\aligned &\int_{D_Y} E_i({az+b\over d}, s)E_i(z,\bar s)dza\\
&=\sum_{j=1}^{\nu(N)} \{ {\delta_{ij}(a/d)^sY^{s+\bar s-1}-
\varphi_{ij}(s)\varphi_{ij}(\bar s)(a/d)^{1-s}Y^{1-s-\bar s}
\over s+\bar s -1}\\
&+\delta_{ij}{(a/d)^s \varphi_{ij}(\bar s) Y^{s-\bar s}
-\varphi_{ij}(s)(a/d)^{1-s}Y^{\bar s-s}\over s-\bar s}\}+o(1)
 \endaligned$$
for nonreal $s$ with $\Re s>1/2$, where $o(1)$ tends to zero
as $Y\to\infty$.  Hence we have
$$\aligned &\int_{D_Y} \left(\sum_{ad=n, 0\leqslant b<d}
E_i({az+b\over d}, s)\right)E_i(z,\bar s)dz\\
&=\{ {Y^{s+\bar s-1}-\sum_{j=1}^{\nu(N)}\varphi_{ij}(s)
\varphi_{ij}(\bar s)Y^{1-s-\bar s}\over s+\bar s -1}
+{\varphi_{ii}(\bar s) Y^{s-\bar s}
-\varphi_{ii}(s)Y^{\bar s-s}\over s-\bar s}\}\\
&\,\,\,\,\times \{\sum_{ad=n, d>0}a^sd^{1-s}\}+o(1).
 \endaligned$$
By partial integration, we obtain
$$h(r)=\frac 1{r^4}\int_0^\infty g^{(4)}(\ln u) u^{ir-1}du $$
for nonzero $r$.   Then it follows that
$$\aligned &\lim_{S\to {1\over 2}^+}\int_{-\infty}^\infty
h(r)\{\int_{D_Y} \left(\sum_{ad=n, 0\leqslant b<d}
E_i({az+b\over d}, S+ir)\right)E_i(z,S-ir)dz\}dr\\
&=\sqrt n\sum_{ad=n, d>0}\{4\pi g(\ln{a\over d})\ln Y
-\sum_{j=1}^{\nu(N)} \int_{-\infty}^\infty h(r)
\left({a\over d}\right)^{ir}\varphi_{ij}^\prime({1\over 2}+ir)
\varphi_{ij}({1\over 2}-ir)dr\\
&+\int_{-\infty}^\infty h(r) \left({a\over d}\right)^{ir}{\varphi_{ii}
({1\over 2}-ir) Y^{2ir}\over ir}dr\}+o(1).\endaligned\tag 3.10$$
By the Riemann-Lebesgue theorem (cf. \S1.8 of \cite{16}), we have
$$\aligned &\lim_{Y\to\infty}\int_{-\infty}^\infty h(r)
\left({a\over d}\right)^{ir}{\varphi_{ii}
({1\over 2}-ir) Y^{2ir}\over ir}dr\\
&=\lim_{Y\to\infty} \int_{-\infty}^\infty h(r)
\Re\left(\varphi_{ii}({1\over 2}-ir)\right)
{\sin\{r\ln(aY^2/d)\}\over r}dr=\pi h(0)
\varphi_{ii}({1\over 2}).\endaligned \tag 3.11$$
By (3.10) and (3.11), we have
$$\aligned \frac 1{\sqrt n}\int_{D_Y}H(z, z)dz
&= \sum_{i=1}^{\nu(N)}\sum_{ad=n, d>0}\{g(\ln{a\over d})\ln Y
+{1\over 4}h(0)\varphi_{ii}({1\over 2})\\
&-\sum_{j=1}^{\nu(N)}{1\over 4\pi}\int_{-\infty}^\infty h(r)
\left({a\over d}\right)^{ir}\varphi_{ij}^\prime({1\over 2}+ir)
\varphi_{ij}({1\over 2}-ir)dr\}+o(1).\endaligned$$
 The stated identity then follows.  \qed\enddemo

  Let $N$ be a square free positive integer.
 It follows from Theorem 3.9 and Theorem 3.10 that
$$\aligned &\lim_{Y\to\infty}\left( c(\infty)_Y
+\sum_{\{P\}, \Gamma_P=\{1_2\}}c(P)_Y\right)\\
&=\sqrt n\delta_n\nu(N)g(0)\ln{\sqrt n\over 2}+{\sqrt n\over 4}
h(0)\left(\delta_n\nu(N)+d(n)\sum_{i=1}^{\nu(N)}
\varphi_{ii}({1\over 2})\right)\\
&+{\sqrt n\over 2}\sum_{w|N, w>0}\sum_{ad=n, d>0, a\neq d}
\sum_b{\ln \{(a-d)^2wN[{C^2\over\ell^2},N]/C^2\}
\over |a-d|}g(\ln \frac ad)\\
&-{\delta_n\nu(N)\sqrt n\over 2\pi}\int_{-\infty}^\infty
h(r)\frac{\Gamma^\prime}\Gamma(1+ir)dr \\
&+\sum_{i,j=1}^{\nu(N)}\sum_{ad=n, d>0}
{\sqrt n\over 4\pi}\int_{-\infty}^\infty h(r)
\left({a\over d}\right)^{ir}\varphi_{ij}^\prime({1\over 2}+ir)
\varphi_{ij}({1\over 2}-ir)dr\\
&+\sum_{ad=n, d>0,a\neq d}{1\over 2}\nu(N)
 |a-d|\int_1^\infty  k\left(\frac{(a-d)^2}nt\right)
\frac{\ln t}{\sqrt {t-1}}dt \endaligned \tag 3.12$$
with $C=(a-d)w-bwN$, where $\ell=(C, bN)$ and the summation
on $b$ is taken over all numbers $b$ such that $C, bN/w\in \Bbb Z, N|C$
 and $0\leqslant b<|a-d|$.   Note that there are exactly $|a-d|$
number of such numbers $b$ by Lemma 3.5.

  Denote by $c(\infty)$ the right side of the identity (3.12).
We conclude that the trace formula (2.3) can be written as
$$\aligned & d(n)h(-\frac i2)+\sqrt n \sum_{j=1}^\infty h(\kappa_j)
\text{tr}_jT_n\\
&= c(I)+\sum_{\{R\}} c(R)+\sum_{\{P\},\,\Gamma_P\neq \{1_2\}}
c(P) +c(\infty) \endaligned \tag 3.13$$
for $\Re s>1$, where the summations on the right side of
the identity are taken over the conjugacy classes.

    \proclaim{Lemma 3.11 (Hejhal \cite{3})}
(cf. Proposition 13.6 of Iwaniec \cite{6})
 Let $N$ be a square free positive integer.
Then, for any pair of cusps $\frak a_i=1/w_i, \frak a_j=1/w_j$
with $N=v_iw_i=v_jw_j$, we have
$$\varphi_{ij}(s)=\varphi(s)p_{ij}(s)$$
where
$$\varphi(s)=\sqrt\pi{\Gamma(s-{1\over 2})\over\Gamma(s)}
{\zeta(2s-1)\over \zeta(2s)}$$
and
$$p_{ij}(s)=\varphi\left((v_i, v_j)(w_i, w_j)\right)
\prod_{p|N}(p^{2s}-1)^{-1}\prod_{p|(w_i, v_j)(w_j, v_i)}
(p^s-p^{1-s}).$$
\endproclaim

   \proclaim{Theorem 3.12}  Let $N$ be square free, and put
$c(P)=\int_{\Gamma_P\backslash \Cal H} k(z, Pz)dz$
for hyperbolic elements $P\in\Gamma^*$.  Then the series
$$\sum_{\{P\}, \, \Gamma_P\neq \{1_2\}} c(P)$$
represents an analytic function in the half-plane $\Re
s> 0$ except for a possible pole at $s=1/2$ and for
possible simple poles at $s=1, \frac 12\pm i\kappa_j$,
$j=1,2,\cdots$. \endproclaim

  \demo{Proof}   We have
$$g^{(4)}(\log u)=A(s)u^{\frac 12-s}+O_s(u^{-\frac 12}),$$
where $A(s)$ is an analytic function of $s$
for $\Re s>0$ and where $O_s(u^{-\frac 12})$ means that, for
every complex number $s$ with
$\Re s>0$, there exists a finite constant $B(s)$ depending
only on $s$ such that
$$|O_s(u^{-\frac 12})|\leqslant B(s)u^{-\frac 12}.$$
Moreover, for every fixed value of $u$, the term
$O_s(u^{-\frac 12})$ also represents an analytic
function of $s$ for $\Re s>0$.  Since
$$h(r)=\frac 1{r^4}\int_0^\infty g^{(4)}(\ln u) u^{ir-1}du $$
for nonzero $r$, we have
$$\aligned h(r)&=\frac 1{r^4}\int_1^\infty g^{(4)}
(\ln u) (u^{ir}+u^{-ir}) \frac {du}u  \\
&=\frac {A(s)}{r^4}\int_1^\infty u^{-\frac 12-s}
(u^{ir}+u^{-ir})du+O_s\left(\frac 1{r^4}\int_1^\infty
u^{-1-\epsilon}du\right)\\
&=\frac {A(s)}{r^4}\left(\frac 1{s-\frac 12-ir}+
\frac 1{s-\frac 12+ir}\right)+O_s(r^{-4})\endaligned$$
for $\Re s>1$ and for nonzero $r$ with $|\Im r|<{1\over 2}-\epsilon$.
By analytic continuation, we obtain that
$$h(r)=\frac{2A(s)(s-\frac 12)}{r^4[(s-\frac 12)^2+
r^2]}+O_s(r^{-4}) \tag 3.14$$
for $\Re s>0$ and for nonzero $r$ with $|\Im r|<{1\over 2}-\epsilon$.
It follows from results of \cite {17}
that the left side of (3.13) is an analytic function of $s$
for $\Re s>0$ except for simple poles at $s=1, \frac 12\pm i\kappa_j$,
$j=1,2,\cdots$.   Then the right
side of (3.13) can be interpreted as an analytic function of $s$ in
the same region by analytic continuation.

   Since $k(t)=(1+t/4)^{-s}$, by Lemma 3.1 we have that $c(R)$ is
analytic for $\Re s>0$ except for a simple pole at $s=1/2$.
There are only a finite number of elliptic conjugacy classes
$\{R\}$.  The term $c(I)$ is a constant.

  Since $g(0)=2\sqrt\pi\Gamma(s-{1\over 2})\Gamma(s)^{-1}$,
$$h(0)=2\sqrt\pi 4^s{\Gamma(s-{1\over 2})\over\Gamma(s)}
\int_1^\infty (u+{1\over u}+2)^{{1\over 2}-s}{du\over u}$$
and
$$g(\ln{a\over d})=2\sqrt\pi 4^{s-{1\over 2}}{\Gamma(s-{1\over 2})
\over \Gamma(s)}\left({a\over d}+{d\over a}+2\right)^{{1\over 2}-s},$$
the sum of first three terms on the right side of the identity (3.12)
is analytic for $\Re s>0$ except for a pole at $s=1/2$.

  Since $k(t)=(1+t/4)^{-s}$,  the sixth term on the right side
of the identity (3.12) is analytic for $\Re s>0$ except for
a pole at $s=1/2$.

 By Stirling's formula the identity
$${\Gamma^\prime(z)\over\Gamma(z)}=\ln z+O(1) \tag 3.15$$
holds uniformly when $|\arg z|\leq \pi-\delta$ for a small positive
number $\delta$.   It follows from (3.14) and (3.15) that the
fourth term on the right side of the identity (3.12)
is analytic for $\Re s>0$ except for a possible pole at $s=1/2$.

  By the functional identity of the Riemann zeta function $\zeta(s)$,
we have $|\varphi(s)|=1$ for $\Re s=1/2$.  This implies that
$$\varphi_{ij}^\prime(s)\varphi_{ij}(\bar s)
={\varphi^\prime(s)\over\varphi(s)}|p_{ij}(s)|^2
+p_{ij}^\prime(s) p_{ij}(\bar s)  \tag 3.16$$
for $\Re s=1/2$ by Lemma 3.11.  The identity
$${\varphi^\prime(s)\over\varphi(s)}=2\ln\pi
-{\Gamma^\prime(s)\over \Gamma(s)}-{\Gamma^\prime(1-s)
\over\Gamma(1-s)}-2{\zeta^\prime(2s)\over\zeta(2s)}
-2{\zeta^\prime(2-2s)\over\zeta(2-2s)} \tag 3.17$$
holds for $\Re s=1/2$.   It follows from (3.14)-(3.17) and
the formula for $p_{ij}(s)$ given in Lemma 3.11 that the fifth term
on the right side of the identity (3.12) is analytic for $\Re s>0$
 except for a possible pole at $s=1/2$.

   Therefore, by (3.13) we have proved that  the series
$$\sum_{\{P\}, \, \Gamma_P\neq \{1_2\}} c(P)$$
represents an analytic function of $s$ in the half-plane
$\Re s>0$ except for a possible pole at $s=1/2$ and for possible simple
poles at $s=1,\frac 12\pm i\kappa_j$, $j=1,2,\cdots$. \qed\enddemo

\heading
4.   Proof of the Main Theorem
\endheading

 \proclaim{Lemma 4.1}   Let $\lambda>0$ be an eigenvalue of
$\Delta$ for $\Gamma_0(N)$ with $N$ square free.  Assume that
$(n, N)=1$.  Put $\tau=1/2+i\kappa$ with $\kappa=\sqrt{\lambda-1/4}$.
Then we have
 $$ 4^\tau\sqrt{\pi n}\frac{\Gamma(i\kappa)}{\Gamma(\tau)}
 \text{tr}T_n=\lim_{s\to\tau}(s-\tau)
\sum_{\{P\}, \,\Gamma_P\neq \{1_2\}}c(P)$$
where the right side is defined as in Theorem 3.12.
\endproclaim

\demo{Proof}  By (2.2), (3.14) and results of \cite {17}, we have
$$\lim_{s\to\tau}(s-\tau)\sum_{j=1, \kappa_j\neq\kappa}^\infty
h(\kappa_j)\text{tr}_jT_n=0.$$
By the proof of Theorem 3.12, we have
$$\lim_{s\to\tau}(s-\tau)\left(d(n)h(-{i\over 2})
-c(I)-c(\infty)-\sum_{\{R\}}c(R)\right)=0.$$
  Then the stated identity follows from
(2.2), (3.13) and Theorem 3.12. \qed\enddemo

A quadratic form $ax^2+bxy+cy^2$, which is denoted by $[a, b, c]$,
is said to be primitive if $(a,b,c)=1$ and $b^2-4ac=d\in\Omega$.
Two quadratic forms $[a, b,c]$ and $[a^\prime, b^\prime,
c^\prime]$ are equivalent if an element $\gamma\in SL_2(\Bbb Z)$
exists such that
$$\left(\matrix {a^\prime}&{b^\prime/2}\\
{b^\prime/2}&{c^\prime}\endmatrix\right)=\gamma^t
\left(\matrix a&{b/2}\\{b/2}&c\endmatrix\right)\gamma,$$
where $\gamma^t$ is the transpose of $\gamma$.  This relation partitions
the quadratic forms into equivalence classes, and two such forms from the
same class have the same discriminant.  The number of classes $h_d$
of primitive indefinite quadratic forms of a given discriminant $d$
is finite, and is called the class number of
indefinite quadratic forms.

  \remark{ Remark}  Siegel \cite {15} proved that
$$\lim_{d\to\infty} \frac{\ln(h_d\ln\epsilon_d)}{\ln d}=\frac 12.\tag 4.1$$
 \endremark

  \proclaim{Lemma 4.2}   $P$ is a hyperbolic element of
$\Gamma^*$ with $\Gamma_P\neq\{1_2\}$ if and only if there
exists a primitive indefinite quadratic form $[a,b,c]$
of discriminant $d$ such that
 $$P=\frac 1{\sqrt n}\left(\matrix
{\frac {v-bNu(a, N)^{-1}}2}&{-cNu(a,N)^{-1}}\\
{aNu(a,N)^{-1}}&{\frac {v+bNu(a,N)^{-1}}2}\endmatrix\right)$$
with
$$v^2-\{dN^2(a,N)^{-2}\}u^2=4n.$$
   If   $\lambda_P$ is an eigenvalue of $P$, then
$$\lambda_P-{1\over\lambda_P}=\pm{Nu\over (a,N)}
{\sqrt d\over \sqrt n }.$$
Let
$$P_0=\left(\matrix {\frac {v_0-bu_0}2}&{-cu_0}\\
{au_0}&{\frac {v_0+bu_0}2}\endmatrix\right),$$
where the pair $(v_0, u_0)$ is the fundamental solution
of Pell's equation $v^2-du^2=4$.   Then $P$ is
$\Gamma_0(N)$-conjugate to a hyperbolic element
$P^\prime\in\Gamma^*$ with $\Gamma_{P^\prime}\neq \{1_2\}$
 if and only if $P_0$ is $\Gamma_0(N)$-conjugate to $P_0^\prime$,
where $P_0^\prime$ is associated with $P^\prime$ similarly as
$P_0$ is associated with $P$. \endproclaim

 \demo{Proof}   Let
$$P=\frac 1{\sqrt n}\left(\matrix A&B\\C&D\endmatrix\right)$$
be a hyperbolic element of $\Gamma^*$ such that $\Gamma_P\neq \{1_2\}$.
Then fixed points $r_1$, $r_2$ of $P$ are not rational numbers
by Lemma 3.3, which satisfy the equation $Cr^2+(D-A)r-B=0$.  This
implies that $\Gamma_P$ is the subgroup of elements in $\Gamma_0(N)$
having $r_1$, $r_2$ as fixed points.
  Let $a=C/\mu, b=(D-A)/\mu$ and $c=-B/\mu$, where $\mu=(C, D-A, -B)$.
Then $[a, b, c]$ is a primitive quadratic form with
 $r_1$, $r_2$ being the roots of the equation $ar^2+br+c=0$.
By Sarnak \cite {9}, the subgroup of elements in $SL_2(\Bbb Z)$
 having $r_1$, $r_2$ as fixed points consists of matrices of the form
$$\left(\matrix {\frac {v-bu}2}&{-cu}\\
{au}&{\frac {v+bu}2}\endmatrix\right)$$
with $v^2-du^2=4$ where $d=b^2-4ac$,  and it is generated by the
primitive hyperbolic element
$$P_0=\left(\matrix {\frac {v_0-bu_0}2}&{-cu_0}\\
{au_0}&{\frac {v_0+bu_0}2}\endmatrix\right)$$
where the pair $(v_0, u_0)$ is the fundamental solution of Pell's equation
$v^2-du^2=4$.

  Since $P$ and $P_0$ have the same fixed points, we have
$A=D-bC/a$ and $B=-cC/a$.  Since $P$ belongs to
$\Gamma^*$ and $AD-BC=n$, $C$ satisfies
$$\cases aD^2-bDC+cC^2=na\\ a|C,\,\, N|C.\endcases\tag 4.2$$
Let $\lambda_P$ be an eigenvalue of $P$.  Then it is a solution of
the equation $\lambda^2-{A+D\over\sqrt n}\lambda+1=0$.
By using $A=D-bC/a$ and $B=-cC/a$, we obtain that
$$\lambda_P-\frac 1{\lambda_P}=\pm
\frac {C\sqrt d}{a\sqrt n }  \tag 4.3$$
and
$$\lambda_P+\frac 1{\lambda_P}=\frac 1{\sqrt n}(2D-\frac ba C). \tag 4.4$$

  Conversely, let a pair $(C, D)$ be a solution of the equation (4.2).
Define $A=D-bC/a$ and $B=-cC/a$.  Then the matrix
$$P=\frac 1{\sqrt n}\left(\matrix A&B\\C&D\endmatrix\right)$$
has the same fixed points as $P_0$, and eigenvalues of
$P$ satisfies (4.3) and (4.4).  We have the decomposition
$$P=\frac 1{\sqrt n}\left(\matrix {\frac n{(D, C)}}&{-\frac ca\alpha C
-(D-\frac ba C)\beta}\\0&{(D, C)}\endmatrix\right)T\cdot T^{-1}
\left(\matrix \alpha &\beta\\{C/(D, C)}&{D/(D, C)}
\endmatrix\right),$$
where $\alpha$ and $\beta$ are integers such that
$\alpha D-\beta C=(D, C)$ and where
$T=\left(\smallmatrix 1&*\\0&1\endsmallmatrix\right)
\in\Gamma_0(N)$ is chosen so that
$$\left(\matrix {n/(D, C)}&*\\0&{(D, C)}\endmatrix\right)T
=\left(\matrix {n/(D, C)}&{-s}\\0&{(D, C)}\endmatrix\right)$$
with $0\leq s<n/(D, C)$.  The first equation of (4.2) can be
written as
$$D^2-bD{C\over a}+cC{C\over a}=n.$$
  Since $(n, N)=1$ and $a|C$, we have $\left(N, (D,C)\right)=1$.
Hence $\left(\smallmatrix \alpha &\beta\\{C/(D, C)}&{D/(D, C)}
\endsmallmatrix\right)$ is an element of $\Gamma_0(N)$.
Therefore, $P$ is an hyperbolic element of $\Gamma^*$
with $\Gamma_P\neq\{1_2\}$ by Lemma 3.3.

  Next, let $v=2D-b\frac Ca$ and $u=\frac Ca$.
Then the equation (4.2) becomes $v^2-du^2=4n $ with $N|au$.
Since $N|au$, this equation can be written as
$$v^2-{dN^2\over (a,N)^2}u^2=4n. \tag 4.5$$
Moreover, we have
$$\frac 1{\sqrt n}\left(\matrix A&B\\C&D\endmatrix\right)
=\frac 1{\sqrt n}\left(\matrix
{\frac {v-bNu(a, N)^{-1}}2}&{-cNu(a,N)^{-1}}\\
{aNu(a,N)^{-1}}&{\frac {v+bNu(a,N)^{-1}}2}
\endmatrix\right). \tag 4.6$$

   Let $P_0^\prime$ be the primitive hyperbolic element of
$SL_2(\Bbb Z)$ corresponding to $[a^\prime, b^\prime, c^\prime]$.
Since the identity
$$\left(\matrix {v_0/2}&0\\ 0&{v_0/2}\endmatrix\right)+
u_0\left(\matrix 0&{-1}\\ 1&0\endmatrix\right)\gamma^t
\left(\matrix a&{b/2}\\{b/2}&c\endmatrix\right)\gamma
=\gamma^{-1}P_0 \gamma$$
holds for every element $\gamma\in\Gamma_0(N)$,
two forms $[a, b, c]$ and $[a^\prime, b^\prime, c^\prime]$
of the same discriminant are equivalent in $\Gamma_0(N)$
if and only if an element $\gamma\in\Gamma_0(N)$ exists
such that $\gamma^{-1}P_0\gamma=P_0^\prime$. \qed\enddemo

  \proclaim{Lemma 4.3}  Let $[a,b,c]$ be a primitive
indefinite quadratic form of discriminant $d$, and let
 $$P=\frac 1{\sqrt n}\left(\matrix
{\frac {v-bNu(a, N)^{-1}}2}&{-cNu(a,N)^{-1}}\\
{aNu(a,N)^{-1}}&{\frac {v+bNu(a,N)^{-1}}2}\endmatrix\right)$$
with $v^2-\{dN^2(a,N)^{-2}\}u^2=4n$
be a hyperbolic element of $\Gamma^*$
with $\Gamma_P\neq \{1_2\}$.  Let $(v_1, u_1)$
with $v_1, u_1>0$ be the
fundamental solution of Pell's equation
$v^2-d_1u^2=4$, where $d_1=dN^2/(a, N)^2$.
Then $\Gamma_P$ is generated by the hyperbolic element
$$P_1=\left(\matrix {{1\over 2}\left(v_1-b{Nu_1\over (a, N)}
\right)}&{-c{Nu_1\over (a, N)}}\\
{a{Nu_1\over (a, N)}}&{{1\over 2}\left(v_1+b{Nu_1\over
(a, N)} \right)}\endmatrix\right)$$
of $\Gamma_0(N)$.\endproclaim

  \demo{Proof}   Let $r_1$ and $r_2$ be the fixed points of $P$.
Then they satisfy $ar^2+br+c=0$, and
$\Gamma_P$ is the subgroup of elements in $\Gamma_0(N)$
having $r_1$, $r_2$ as fixed points.  By Sarnak \cite {9}, the
subgroup of elements in $SL_2(\Bbb Z)$ having $r_1$, $r_2$ as
fixed points consists of matrices of the form
$$\left(\matrix {\frac {v-bu}2}&{-cu}\\
{au}&{\frac {v+bu}2}\endmatrix\right)$$
with $v^2-du^2=4$ where $d=b^2-4ac$,  and it is generated
by the primitive hyperbolic element
$$P_0=\left(\matrix {\frac {v_0-bu_0}2}&{-cu_0}\\
{au_0}&{\frac {v_0+bu_0}2}\endmatrix\right)$$
where the pair $(v_0, u_0)$ is the fundamental solution
of Pell's equation $v^2-du^2=4$.

  Since $\Gamma_P$ is cyclic, a solution $v_1, u_1^\prime>0$
of Pell's equation $v^2-du^2=4$ exists such that
 $\Gamma_P$ is generated by
$$P_1=\left(\matrix {\frac {v_1-bu_1^\prime}2}&{-cu_1^\prime}\\
{au_1^\prime}&{\frac {v_1+bu_1^\prime}2}\endmatrix\right)$$
and such that $P_1$ is the smallest positive integer
power of $P_0$ among all powers of $P_0$ belonging to
$\Gamma_0(N)$.   Note that the eigenvalues of $P_1$ are
$${v_1\pm \sqrt d u_1^\prime\over 2}.$$
Since $\Gamma_P$ is generated by $P_1$, $(v_1, u_1^\prime)$
is the minimal solution of the equation
$v^2-du^2=4$ with $N|au_1^\prime$ in the sense that
$(v_1+\sqrt du_1^\prime)/2$ is of the smallest value
among all such solutions.  Since $N|au_1^\prime$,
we have ${N\over (a, N)}|u_1^\prime$.  Write
$$u_1^\prime={Nu_1\over (a, N)}.$$
Then the pair $(v_1, u_1)$ with $v_1, u_1>0$ must
be the fundamental solution of the Pell equation
$v^2-d_1u^2=4$, where $d_1=dN^2(a, N)^{-2}$.
The stated result then follows. \qed\enddemo

 Two quadratic forms $[a, b,c]$ and $[a^\prime, b^\prime,
c^\prime]$ are equivalent in $\Gamma_0(N)$ if an element
$\gamma\in \Gamma_0(N)$ exists such that
$$\left(\matrix {a^\prime}&{b^\prime/2}\\
{b^\prime/2}&{c^\prime}\endmatrix\right)=\gamma^t
\left(\matrix a&{b/2}\\{b/2}&c\endmatrix\right)\gamma.$$
This relation partitions the quadratic forms into
equivalence classes, and two such forms from the
same class have the same discriminant.  The number of such classes of a
given discriminant $d$ is finite, and is denoted by $H_d$.

  \proclaim{Lemma 4.4}   Let $[a_j,b_j,c_j]$, $j=1,2,\cdots,H_d$,
be a set of representatives for classes of primitive indefinite
quadratic forms of discriminant $d$, which are not equivalent
under $\Gamma_0(N)$.  Then we have
$$\aligned &\sum_{\{P\}:\,\,P\in\Gamma^*,\,
\Gamma_P\neq \{1_2\}}c(P)\\
&= 4\sqrt{\pi n}{\Gamma(s-{1\over 2})\over N\Gamma(s)}
\sum_{d\in \Omega}\sum_{j=1}^{H_d}\sum_u
{(a_j,N)\ln \epsilon_{d_1}\over u\sqrt d}
\left(1+{d(Nu)^2\over 4n(a_j,N)^2}\right)^{{1\over 2}-s}
\endaligned$$
for $\Re s>1$, where $d_1=dN^2(a_j, N)^{-2}$ and
where the summation on $u$ is taken over all the
positive integers $u$ such that $4n+dN^2(a_j,N)^{-2}u^2$ is
the square of an integer.  \endproclaim

 \demo{Proof}  It follows from (2.1) and Theorem 3.2 that
$$\sum_{\{P\}:\,\,P\in\Gamma^*,\, \Gamma_P\neq \{1_2\}}c(P)
=2\sqrt\pi{\Gamma(s-{1\over 2})\over\Gamma(s)}\sum_{\{P\}}
{\ln NP_1\over \lambda_P-1/\lambda_P}\left(1+{(\lambda_P
-1/\lambda_P)^2\over 4}\right)^{{1\over 2}-s}$$
for $\Re s>1$, where $\lambda_P>1$ is an eigenvalue of $P$
and $P_1$ is given in Lemma 4.3.  Let $P$ be associated with
a primitive indefinite quadratic form $[a_j, b_j,c_j]$
as in Lemma 4.2.  Then by Lemma 4.2, we have
$$\lambda_P-{1\over\lambda_P}={Nu\over (a_j,N)}
{\sqrt d\over \sqrt n }.$$
By Lemma 4.3, we have
$$\sqrt{NP_1}={v_1+\sqrt{d_1}u_1\over 2}=\epsilon_{d_1}.$$

  If $P^\prime$ is a hyperbolic of $\Gamma^*$ with
$\Gamma_{P^\prime}\neq\{1_2\}$, and is associated with
a primitive indefinite quadratic form
$[a_{j^\prime}, b_{j^\prime},c_{j^\prime}]$
as in Lemma 4.2, then $P$ and $P^\prime$ are
$\Gamma_0(N)$-conjugate if and only if
$[a_j, b_j,c_j]$ and $[a_{j^\prime}, b_{j^\prime},c_{j^\prime}]$
are $\Gamma_0(N)$-conjugate by the last statement of
Lemma 4.2.  Next, let $T$ be a hyperbolic of $\Gamma^*$ with
$\Gamma_T\neq\{1_2\}$.  Assume that $T$ is associated with
a primitive indefinite quadratic form $[a,b,c]$
as in Lemma 4.2.  If the discriminant of $[a,b,c]$ is
not equal to the discriminant of $[a_j, b_j,c_j]$
which is associated with $P$, then $P$ and $T$ are
not $\Gamma_0(N)$-conjugate.  The stated identity
then follows.  \qed\enddemo

  \proclaim{Lemma 4.5}  Let $N$ be square free, and
let $k$ be a divisor of $N$.  Then
the number of indefinite primitive quadratic forms
$[a,b,c]$ with $(a,N)=k$ of discriminant $d$, which are not
equivalent under $\Gamma_0(N)$, is equal to
$$h_{d_1}\prod_{p|k}\left(1+\left({d\over p}\right)\right)$$
where $d_1=dN^2/k^2$.  \endproclaim

  \demo{Proof}   Let $[a,b,c]$ and $[a^\prime, b^\prime, c^\prime]$
be two indefinite primitive quadratic forms of discriminant $d$
with $(a,N)=k=(a^\prime, N)$.  If they are equivalent under
$\Gamma_0(N)$, then an element $\left(\smallmatrix\alpha&\beta
\\ \gamma&\delta\endsmallmatrix\right)\in \Gamma_0(N)$ exists
such that
$${b^\prime\over 2}=\alpha(a\beta+{b\over 2}\delta)
+\gamma({b\over 2}\beta+c\delta).$$
This implies that $b^\prime\equiv b$ (mod $2k$).
In particular, if $b^\prime\not\equiv b$ (mod $2k$),
then $[a,b,c]$ and $[a^\prime, b^\prime, c^\prime]$
are not $\Gamma_0(N)$-equivalent.

  Assume that $\varrho$ is an integer with
$1\leqslant \varrho\leqslant 2k$.
 Denote by $\Lambda_{k,d,\varrho}$ the set of
representatives of indefinite primitive quadratic forms
$[a,b,c]$ of discriminant $d$ with $(a,N)=k$ and
$b\equiv\varrho$ (mod $2k$), which are not equivalent under
$\Gamma_0(N)$.    Let $\Cal L_{N,d_1,\varrho N/k}^0$ be
the set of representatives of indefinite quadratic forms
$[aN,b,c]$ of discriminant $d_1$ with $(a,b,c)=1, (N,b,c)=N/k$
and $b\equiv \varrho N/k$ (mod $2N$),
which are not equivalent under $\Gamma_0(N)$.

   A map from $\Lambda_{k, d,\varrho}$ to
$\Cal L_{N,d_1,\varrho N/k}^0$ is defined by
$T: [a,b,c]\to [aN/k, bN/k, cN/k]$.  We claim that $T$ is
bijective.  By the definition of $\Lambda_{k,d,\varrho}$,
we see that $T$ is injective.   Conversely, if $[a_1N,b_1,c_1]$
is an element of $\Cal L_{N,d_1,\varrho N/k}^0$, then we have
$(a_1,b_1, c_1)=1, (N, b_1, c_1)=N/k,
 b_1\equiv \varrho N/k$ (mod $2N$) and
$b_1^2-4Na_1c_1=d_1$.  Let $[a_0,b_0,c_0]$ be an element of
$\Lambda_{k, d,\varrho}$.  Then
$$b_1^2=(b_0N/k)^2+4\ell b_0 N{N\over k}
+4\ell^2 N^2$$
for some integer $\ell$.  Since
$$d_1=(b_0N/k)^2-4{a_0\over k}N\cdot c_0{N\over k},$$
we have
$$(b_0N/k)^2-4{a_0\over k}N\cdot c_0{N\over k}=d_1
=(b_0N/k)^2+4\ell b_0 N{N\over k}+4\ell^2 N^2-4Na_1c_1.$$
That is, the identity
$$-{a_0\over k}\cdot c_0{N\over k}
=\ell b_0 {N\over k}+\ell^2 N-a_1c_1 \tag 4.7$$
holds for some integer $\ell$.  Note that $(a_0, N)=k$.
Since $(N, b_1, c_1)=N/k$ and $(a_1, b_1, c_1)=1$,
we have $(a_1, N/k)=1$, and hence it follows from (4.7) that
$N/k|c_1$.  Let $a=a_1 k, b=b_1 k/N$ and $c=c_1 k/N$.
Then we have $(a,N)=k, b\equiv\varrho$ (mod $2k$) and $d=b^2-4ac$.
We claim that $(a,b,c)=1$, that is, $(a_1k, b_1k/N, c_1 k/N)=1$.
Since $(a_1, b_1, c_1)=1$, it is enough to show that
$(k, b_1, c_1)=1$.   Since $(N, b_1, c_1)=N/k$, we
must have $(k, b_1, c_1)=1$,
and therefore we have $(a,b,c)=1$.  Thus, $[a,b,c]$ is an element
of $\Lambda_{k, d,\varrho}$,  and $T$ maps it into
the element $[a_1N,b_1,c_1]$ of $\Cal L_{N,d_1,\varrho N/k}^0$.
Therefore, $T$ is surjective.   Thus, we have proved that
$T$ is a bijection if the set $\Lambda_{k, d,\varrho}$ is
not empty.  By Proposition, p.\,505 of Gross, Kohnen
and Zagier \cite{2}, the number of elements contained
in $\Cal L_{N,d_1,\varrho N/k}^0$
is $h_{d_1}$, and hence the set $\Lambda_{k, d,\varrho}$
contains $h_{d_1}$ elements if it is not empty.

If $\varrho_1\not\equiv\varrho_2$
(mod\,$2k$), then the set of $\Gamma_0(N)$-equivalence
classes represented by elements in $\Lambda_{k,d,\varrho_1}$
is disjoint from the set of $\Gamma_0(N)$-equivalence
classes represented by elements in $\Lambda_{k,d,\varrho_2}$.
Now, we want to count the number of non-empty sets
$\Lambda_{k, d,\varrho}$.  That is, we want to count the
number of solutions $\varrho$ of the equation
$$\varrho^2\equiv d\,\, (\text{mod}\, 4k),\,\,
1\leqslant \varrho\leqslant 2k. \tag 4.8$$

If $(d, k)=1$, by Theorem 3.4 of Chapter 12, Hua \cite{4}
the number of solution of the equation (4.8) is equal to
$$\prod_{p|k}\left(1+\left({d\over p}\right)\right).$$
Next, we consider the case when there exists a prime number
$q$ satisfying $q|k$ and $q^2|d$.  Then we have
$q|\varrho$, and the equation (4.8) can be written as
 $$\left({\varrho\over q}\right)^2\equiv {d\over q^2}
\,\, (\text{mod}\, {4k\over q}),\,\,
1\leqslant {\varrho\over q}\leqslant {2k\over q}. \tag 4.9$$
Dividing out $\varrho, d$ and $k$ by all such prime numbers
$q$ as in (4.9), we can reduce the second case to the first
case when $(d, k)=1$.  Then, by using properties of the
Legendre symbol we obtain that the number of solution of
the equation (4.8) in the second case is still equal to
$$\prod_{p|k}\left(1+\left({d\over p}\right)\right).$$
Finally, we consider the case when there exists a prime number
$q$ satisfying $q|k, q|d$ and $q^2\nmid d$.  We have again
$q|\varrho$.  The equation (4.8) can be written as
 $$q\left({\varrho\over q}\right)^2\equiv {d\over q}
\,\, (\text{mod}\, {4k\over q}),\,\,
1\leqslant {\varrho\over q}\leqslant {2k\over q}. \tag 4.10$$
We can assume that $q\neq 2$.  Otherwise, if $q=2$ then we
must have $q^2|d$.  Then we have $(q, 4k/q)=1$, and hence a
number $x_q$ exists such that $qx_q\equiv 1$ (mod $4k/q$).
   Then the equation (4.10) can be written as
$$\left({\varrho\over q}\right)^2\equiv x_q{d\over q}
\,\, (\text{mod}\, {4k\over q}),\,\,
1\leqslant {\varrho\over q}\leqslant {2k\over q}. \tag 4.11$$
Note that we have
$$\left({x_q d/q\over p}\right)=\left({q\over p}\right)
\left({d/q\over p}\right)=\left({d\over p}\right) \tag 4.12$$
for any prime number $p|(k/q)$.  Dividing out all such
primes $q$ as in (4.11) and using (4.12), we obtain that
the number of solution of the equation (4.8) in the final
 case is equal to
$$\prod_{p|k}\left(1+\left({d\over p}\right)\right).$$
This completes the proof of the lemma.  \qed\enddemo

 By Lemma 4.5, we get the following corollary.

  \proclaim{Corollary}   Let $H_d$ be the number of
indefinite primitive quadratic forms of discriminant $d$,
which are not equivalent under $\Gamma_0(N)$.  Then we have
 $$H_d=\sum_{k|N}h_{d(N/k)^2}\prod_{p|k}
\left(1+\left({d\over p}\right)\right).$$ \endproclaim

  \proclaim{Theorem 4.6}   Let $N$ be a square free positive integer
with $(n, N)=1$.   Then we have
$$\aligned &\sum_{\{P\}:\,\,P\in\Gamma^*,\,
\Gamma_P\neq \{1_2\}}c(P)
\\&= 4\sqrt{\pi n}{\Gamma(s-{1\over 2})
\over \Gamma(s)}
 \sum_{m|N}\sum_{k|N}{\mu((m, k))\over (m,k)}
\sum_{d\in \Omega}\sum_u \left({d\over m}\right)
{h_d\ln \epsilon_d\over u\sqrt d}
\left(1+{d k^2u^2\over 4n}\right)^{{1\over 2}-s}\endaligned $$
for $\Re s>1$, where the summation on $u$ is taken over
all the positive integers $u$ such that
$\sqrt{4n+d k^2 u^2}\in\Bbb Z$.\endproclaim

  \demo{Proof}  Let $[a_j,b_j,c_j]$, $j=1,2,\cdots,H_d$,
be a set of representatives for classes of primitive indefinite
quadratic forms of discriminant $d$, which are not equivalent
under $\Gamma_0(N)$.  By Lemma 4.4, we have
$$\sum_{\{P\}:\,\,P\in\Gamma^*,\,
\Gamma_P\neq \{1_2\}}c(P)
= 4\sqrt{\pi n}{\Gamma(s-{1\over 2})\over \Gamma(s)}
\sum_{d\in \Omega}\sum_{j=1}^{H_d}\sum_u
{\ln \epsilon_{d_1}\over u\sqrt {d_1}}
\left(1+{d_1 u^2\over 4n}\right)^{{1\over 2}-s}$$
for $\Re s>1$, where $d_1=dN^2(a_j, N)^{-2}$ and
where the summation on $u$ is taken over all the
positive integers $u$ such that $\sqrt{4n+d_1 u^2}
\in\Bbb Z$.   Denote $k=N/(a_j, N)$.  Then $k|N$.
By using Lemma 4.5, we can write the above identity as
$$\aligned &\sum_{\{P\}:\,\,P\in\Gamma^*,\,
\Gamma_P\neq \{1_2\}}c(P)\\
&= 4\sqrt{\pi n}{\Gamma(s-{1\over 2})\over \Gamma(s)}
\sum_{k|N}\sum_{d\in \Omega}\sum_u
\prod_{p|{N\over k}}\left(1+\left({d\over p}\right)
\right){h_{d_1}\ln \epsilon_{d_1}\over u\sqrt {d_1}}
\left(1+{d_1u^2\over 4n}\right)^{{1\over 2}-s}
\endaligned \tag 4.13$$
for $\Re s>1$, where $d_1=dk^2$.   By using Dirichlet's
class number formula
$$h_{d_1}\ln\epsilon_{d_1}=\sqrt{d_1}L(1, \chi_{d_1})$$
and by using the identity (See Theorem 11.2 of Chapter 12,
Hua \cite{4})
$$L(1, \chi_{d_1})=L(1, \chi_d)\prod_{p|k}
\left(1-\left({d\over p}\right)p^{-1}\right),$$
we can write (4.13) as
$$\aligned &\sum_{\{P\}:\,\,P\in\Gamma^*,\,
\Gamma_P\neq \{1_2\}}c(P)= 4\sqrt{\pi n}{\Gamma(s-{1\over 2})
\over \Gamma(s)}\\
&\times \sum_{k|N}\sum_{d\in \Omega}\sum_u
\prod_{p|{N\over k}}\left(1+\left({d\over p}\right)
\right)\prod_{p|k}\left(1-\left({d\over p}\right)p^{-1}\right)
{h_d\ln \epsilon_d\over u\sqrt d}
\left(1+{d k^2u^2\over 4n}\right)^{{1\over 2}-s}\endaligned $$
for $\Re s>1$, where the summation on $u$ is taken over
all the positive integers $u$ such that
$\sqrt{4n+d k^2 u^2}\in\Bbb Z$.  Since
$$\prod_{p|{N\over k}}\left(1+\left({d\over p}\right)
\right)\prod_{p|k}\left(1-\left({d\over p}\right)p^{-1}\right)
=\sum_{m|N} \left({d\over m}\right){\mu((m,k))\over (m,k)},$$
we have
$$\aligned &\sum_{\{P\}:\,\,P\in\Gamma^*,\,
\Gamma_P\neq \{1_2\}}c(P)= 4\sqrt{\pi n}{\Gamma(s-{1\over 2})
\over \Gamma(s)}\\
&\times \sum_{m|N}\sum_{k|N}{\mu((m,k))\over (m,k)}
\sum_{d\in \Omega}\sum_u \left({d\over m}\right)
{h_d\ln \epsilon_d\over u\sqrt d}
\left(1+{d k^2u^2\over 4n}\right)^{{1\over 2}-s}.\endaligned
\tag 4.14$$

  Next, we want to show that
$$\sum_{d\in \Omega}\sum_u \left({d\over m}\right)
{h_d\ln \epsilon_d\over u\sqrt d}
\left(1+{d k^2u^2\over 4n}\right)^{{1\over 2}-s}$$
is absolutely convergent for $\sigma=\Re s>1$.  Since
$$|\sum_{d\in \Omega}\sum_u \left({d\over m}\right)
{h_d\ln \epsilon_d\over u\sqrt d}
\left(1+{d k^2u^2\over 4n}\right)^{{1\over 2}-s}|
\leqslant \sum_{d\in \Omega}\sum_u
{h_d\ln \epsilon_d\over u\sqrt d}
\left(1+{d k^2u^2\over 4n}\right)^{{1\over 2}-\sigma}.\tag 4.15$$
It is proved in Li \cite{8} that the right side of (4.15)
is convergent for $\sigma>1$, and hence, the right side of
the stated identity is absolutely convergent for
$\Re s>1$.

This completes the proof of the theorem. \qed\enddemo

\vskip0.15truein

 \demo{Proof of Theorem 1}  It is proved at the end of
Li \cite{8} that
$$\sum_{d\in\Omega, u} \left|\frac{h_d\ln \epsilon_d} {\sqrt d u}
(1+\frac {du^2}{4n})^{\frac 12-\sigma}-
\frac{h_d\ln \epsilon_d} {(du^2)^\sigma}\right|\ll \sum_{d\in\Omega, u}
\frac{(du^2)^{\frac 12+\epsilon-1-\sigma}}{u^{1+2\epsilon}}<\infty$$
for $\sigma>0$.  Then it follows from Theorem 4.6 that
$$\aligned & \lim_{s\to 1/2+i\kappa}(s-{1\over 2}-i\kappa)\sum_{\{P\}, \,
\Gamma_P\neq \{1_2\}}c(P)=\frac{(4\pi)^{1/2}\Gamma(i\kappa)}
{(4n)^{-1/2-i\kappa}\Gamma(1/2+i\kappa)}\\
&\times \sum_{m|N}\sum_{k|N}k^{-2i\kappa}\prod_{p|k}
(1-p^{-1})\lim_{s\to 1/2+i\kappa}(s-{1\over 2}-i\kappa)\sum_{d\in \Omega}\sum_u
\left({d\over m}\right){h_d\ln \epsilon_d\over (du^2)^s}.
\endaligned \tag 4.16 $$
Theorem 3.12 shows that the function on the right side of (4.16)
represents an analytic function in the half-plane $\Re
s> 0$ except for a possible pole at $s=1/2$ and for
possible simple poles at $s=1, \frac 12\pm i\kappa_j$,
$j=1,2,\cdots$.  The stated identity then follows from Lemma 4.1.

  This completes the proof of the theorem.  \enddemo

\enddocument